\documentclass{amsart}

\usepackage{amsmath,amssymb}
\usepackage{amsthm}
\usepackage{amsrefs}
\usepackage{indentfirst}
\usepackage{mathrsfs}
\usepackage{graphicx}
\usepackage{hyperref}
\usepackage{xcolor}
\usepackage{fourier}
\usepackage[margin=1.5in]{geometry}
\usepackage[ruled,norelsize,linesnumbered,vlined,algo2e]{algorithm2e}

\theoremstyle{plain}

\theoremstyle{definition}

\theoremstyle{remark}

\newcommand{\ud}{\,\mathrm{d}}

\newcommand{\RR}{\mathbb{R}}
\newcommand{\NN}{\mathbb{N}}
\newcommand{\ZZ}{\mathbb{Z}}

\newcommand{\Or}{\mathcal{O}}

\newcommand{\wt}[1]{\widetilde{#1}}
\newcommand{\wh}[1]{\widehat{#1}}

\IfFileExists{mathabx.sty}%
  {\DeclareFontFamily{U}{mathx}{\hyphenchar\font45}%
   \DeclareFontShape{U}{mathx}{m}{n}{<->mathx10}{}%
   \DeclareSymbolFont{mathx}{U}{mathx}{m}{n}%
   \DeclareFontSubstitution{U}{mathx}{m}{n}%
   \DeclareMathAccent{\widebar}{0}{mathx}{"73}%
}{%
  \PackageWarning{mathabx}{%
    Package mathabx not available, therefore\MessageBreak substituting
    widebar with overline\MessageBreak }%
  \newcommand{\widebar}[1]{\overline{#1}}%
}
\newcommand{\wb}[1]{\widebar{#1}}

\newcommand{\mc}[1]{\mathcal{#1}}

\newcommand{\abs}[1]{\lvert#1\rvert}

\newcommand{\norm}[1]{\lVert#1\rVert}

\newcommand{\grid}{\text{grid}}
\newcommand{\col}{\text{col}}

\title{Fast algorithm for periodic density fitting for Bloch waves}

\author{Jianfeng Lu} 
\address{Department of Mathematics, Department of
  Physics, and Department of Chemistry, Duke University, Box 90320,
  Durham NC 27708, USA} 
\email{jianfeng@math.duke.edu}

\author{Lexing Ying}
\address{Department of Mathematics and Institute of Computational and Mathematical Engineering, Stanford University, 450 Serra Mall, Bldg 380, Stanford CA 94305, USA}
\email{lexing@stanford.edu}

\date{\today} \thanks{The work of J.L. is supported in part by the
  Alfred P.~Sloan Foundation and the National Science Foundation under
  grant DMS-1312659 and DMS-1454939. The work of L.Y. is partially
  supported by the National Science Foundation under grant DMS-1521830
  and the U.S. Department of Energy's Advanced Scientific Computing
  Research program under grant DE-FC02-13ER26134/DE-SC0009409. We
  would like to thank helpful discussions with Volker Blum and Xinguo
  Ren.}

\begin{document}

\begin{abstract}
  We propose an efficient algorithm for density fitting of Bloch waves
  for Hamiltonian operators with periodic potential. The algorithm is
  based on column selection and random Fourier projection of the
  orbital functions. The computational cost of the algorithm scales as
  $\Or\bigl(N_{\grid} N^2 + N_{\grid} NK \log (NK)\bigr)$, where
  $N_{\grid}$ is number of spatial grid points, $K$ is the number of
  sampling $k$-points in first Brillouin zone, and $N$ is the number
  of bands under consideration. We validate the algorithm by numerical
  examples in both two and three dimensions.
\end{abstract}
\maketitle

\section{Introduction}

Consider a Hamiltonian operator with periodic potential
\begin{equation}
  H = - \frac{1}{2}\Delta + V. 
\end{equation}
For simplicity of notation, we assume without loss of generality that
$V$ is periodic with respect to $\ZZ^d$, where $d$ is the spatial
dimension: \textit{i.e.}, $V(x + e_i) = V(x)$ for
$\{e_i, \, i = 1, \ldots, d\}$ the set of standard Cartesian basis vectors
of $\RR^d$. From the standard Bloch-Floquet theory (see \textit{e.g.},
\cite{ReedSimon4}), the spectrum of the operator $H$ is given by the
Bloch eigenvalue problem on the unit cell $\Gamma = [0, 1)^d$:
\begin{subequations}
\begin{align}
  & (-\frac{1}{2}\Delta + V) \psi_{n, k} = E_{n, k} \psi_{n, k}, \qquad \text{on } \Gamma; \\
  & e^{-i (k \cdot x)} \psi_{n, k}(x) \quad \text{is periodic on $\Gamma$}. \label{eq:blochperiodic}
\end{align}
\end{subequations}
Here $k \in \Gamma^{\ast} = [-\pi, \pi)^d$ is known as the crystal
momentum, where $\Gamma^{\ast}$ is called the first Brillouin zone
(FBZ), which is the unit cell of the reciprocal lattice.  The
eigenfunction $\psi_{n, k}$ for $n \in \NN$ and $k \in \Gamma^{\ast}$
are called Bloch waves. We extend them from $\Gamma$ to the whole
$\RR^3$ using the Bloch-periodicity \eqref{eq:blochperiodic}.  It
would be convenient to introduce the periodic part of the Bloch
function, denoted as
\begin{equation}
  u_{n, k}(x) = e^{-i (k \cdot x)} \psi_{n, k}(x), \qquad x \in \RR^d. 
\end{equation}
By definition, $u_{n, k}$ is periodic with respect to $\ZZ^d$. 

In this work, we are interested in the pair density $\wb{u}_{n,
  k}(x) u_{m, l}(x)$ for $n, m \in \NN, k, l \in
\Gamma^{\ast}$, and $x \in \Gamma$. In particular, we aim at an approximation of the pair density of the type 
\begin{equation}\label{eq:dfapprox}
  \wb{u}_{n, k}(x) u_{m, l}(x) \approx \sum_{\mu} \wb{C}_{n,k}^{\mu} C_{m, l}^{\mu} P_{\mu}(x), 
\end{equation}
where $\{P_{\mu}\}$ is a set of auxiliary basis functions we use to expand the pair density and $C_{n,k}^{\mu}$ gives the expansion coefficients. 
In terms of the Bloch waves, the approximation is then 
\begin{equation}
  \wb{\psi}_{n, k}(x) \psi_{m, l}(x) \approx \sum_{\mu} \wb{C}_{n,k}^{\mu} C_{m, l}^{\mu} P_{\mu}(x) e^{i (l - k) \cdot x}. 
\end{equation}

The approximation of pair density in the form of 
\begin{equation}
  \wb{u}_{n,k}(x) u_{m,l}(x) \approx \sum_{\mu} T_{nk,ml}^{\mu} P_{\mu}(x)
\end{equation}
is known as periodic density fitting or density fitting for crystals
in the literature \cites{Maschio:07, Lorenz:11, Levchenko:15}. Note
that the approximation \eqref{eq:dfapprox} is a special form, in which
the expansion coefficient
$T_{nk,ml}^{\mu} = \wb{C}_{n,k}^{\mu} C_{m,l}^{\mu}$ has a separable
structure. The periodic density fitting is a generalization of the
conventional density fitting for molecules (see \textit{e.g.},
\cites{Ren_etal:12, DunlapConnollySabin:79, Schutz_etal:10,
  SodtSubotnikHeadGordon:06, Vahtras:93, Weigend:98}), which has wide
applications in electronic structure calculations. In density fitting,
the pair density for a set of orbital functions $\{\psi_n\}$ (without
the $k$-dependence) is approximated as
\begin{equation}
  \wb{\psi}_n(x) \psi_m(x) \approx \sum_{\mu} T_{nm}^{\mu} P_{\mu}(x). 
\end{equation}
Here with slight abuse of notation, we use again $T$ to denote the
expansion coefficients. 

Conventionally in density fitting, a set of auxiliary functions
$P_{\mu}(x)$ are chosen, for example as Gaussian functions or atom
centered orbital functions, the expansion coefficients $T$ are
obtained using least square fitting with respect to a suitable metric
($L^2$ or Coulomb metric, see Section~\ref{sec:numerics} below). The
computational cost to obtain these coefficients scales as $\Or(N^4)$,
where $N$ is the total number of orbital functions. Moreover, the
conventional procedure does not yield the separable structure, which
is useful to reduce the computational cost for electronic structure
calculations, see e.g., \cites{Hohenstein:13,
  HohensteinParrishMartinez:12, ParrishHohensteinMartinez:12,
  Parrish:14, Shenvi:13, Shenvi:14}.

In our previous work \cite{LuYing:15}, we developed a
$\Or(N_{\grid} N^2 \log N)$ scaling algorithm for density fitting
without the $k$ point dependence, where $N_{\grid}$ is number of
spatial grids and $N$ is the number of orbital functions.  The novelty
is to recast the density fitting problem as a column selection
procedure: Instead of fixing the auxiliary basis functions and compute
the expansion coefficients, we look for a collection of grid points
$\{x_{\mu}\}$ that well represent the pair densities, such that the
auxiliary basis functions follow as a least square fitting.
\begin{equation}
  \wb{\psi}_m(x) \psi_n(x) \approx \sum_{\mu} \wb{\psi}_m(x_{\mu}) 
  \psi_n(x_{\mu}) P_{\mu}(x). 
\end{equation}

This work extends our approach to periodic density fitting of Bloch
waves. In order to further reduce the computational cost, we use a
more efficient scheme for random projection based on the tensor
product structure of the pair density $\rho$. Our resulting density
fitting scheme scales as
$\Or\bigl(N_{\grid} N^2 + N_{\grid} NK \log (NK)\bigr)$, where
$N_{\grid}$ is number of spatial grid points, $K$ is the number of
$k$-points in first Brillouin zone, and $N$ is the number of bands
under consideration. Recall that in the periodic case, the total
number of Bloch waves is given by $NK$. Note that if we take $K = 1$
so that we go back to the situation considered in \cite{LuYing:15},
our current algorithm costs as $\Or(N_{\grid} N^2)$ and is hence more
efficient than our previous algorithm.  Moreover, the computational
cost scales almost linearly with respect to the increase of $K$, with
a ``quadratic'' coefficient $\Or(N_{\grid} N)$; in fact, in practice,
the actual running time is dominated by the $\Or(N_{\grid} N^2)$ term, independent of $K$.

\section{Algorithms}\label{sec:algorithm}

Let us consider a set of given orbitals $u_{n,k}(x)$ for $n = 1,
\ldots, N$, $k \in \mc{K}$ a discretization of the first Brillouin
zone with $\abs{\mc{K}} = K$, and $x \in \mc{X}$ a grid in the unit
cell with $\abs{\mc{X}} = N_{\grid}$.  Similar to our previous
work \cite{LuYing:15}, the basic idea is to select a suitable subset
of grid points to represent the pair densities $\rho_{nkml}(x) =
\wb{u}_{n, k}(x) u_{m,l}(x)$.
\begin{equation}
  \rho_{nkml}(x) \approx \sum_{\mu} \rho_{nkml}(x_{\mu}) P_{\mu}(x)
  = \sum_{\mu} \wb{u}_{n,k}(x_{\mu}) u_{m,l}(x_{\mu}) P_{\mu}(x). 
\end{equation}
After the grid points $\{x_{\mu}\}$ are determined, the auxiliary
basis functions $P_{\mu}$ follow from a least square fitting.

Note that if we view $\rho$ as a matrix of dimension $N^2K^2 \times
N_{\grid}$ with $(nkml)$ being the row index and $x$ being the column
index, the choice of $x_{\mu}$ amounts to select a sub-collection of
the columns to represent all the matrix columns.  One immediate
approach for the column selection is to use a pivoted QR algorithm
\cite{GolubVanLoan} on $\rho$. For completeness, we recall the column
selection algorithm based on pivoted QR algorithm here. 

\begin{algorithm2e}[H]
\SetKwInOut{Input}{Input}
\SetKwInOut{Output}{Output}

\Input{$m\times n$ matrix $M$, error tolerance $\mathsf{tol}$}
\Output{An $m \times N_{\col}$ submatrix $\wt{M}$ of $M$ and an
  $N_{\col} \times n$ matrix $P$, such that $M \approx \wt{M} P$}

Compute the pivoted QR decomposition $[Q, R, E] = \mbox{qr}(M)$, so that 
\begin{equation*}
  QR = ME,
\end{equation*}
where $E$ is an $n\times n$ permutation matrix, $Q$ is an $m \times m$
unitary matrix and $R$ is an $m \times n$ upper triangular matrix with
diagonal entries in decreasing order;

Set $N_{\col}$ such that 
\begin{equation*}
  R_{N_{\col}, N_{\col}} \geq \mathsf{tol} \cdot R_{1, 1} 
  > R_{N_{\col}+1, N_{\col}+1};
\end{equation*}

Set $\wt{M} = (ME)_{:, 1:N_{\col}}$, the first $N_{\col}$ columns of
$ME$, where \textsf{Matlab} notations are used for submatrices. Note
that $ME$ amounts to a permutation of the columns of $M$.

Compute $P = R_{1:N_{\col}, 1:N_{\col}}^{-1} R_{1:N_{\col}, :}
E^{-1}$.

\caption{Column selection based on pivoted QR}\label{alg:cs}
\end{algorithm2e}
The computational cost of directly applying Algorithm~\ref{alg:cs} on
the matrix $\rho$ is however prohibitively expensive, since the
computational complexity of the QR step scales as
$\Or(N_{\grid} N^4K^4)$ (recall that the matrix $\rho$ has dimension
$N^2K^2 \times N_{\grid}$).

In order to reduce the complexity, it is thus crucial to apply a
random projection on the matrix $\rho$, in the spirit of the
pioneering works \cites{Liberty:07, Woolfe:08}. The idea is to find an
alternative representation of the column space of $\rho$ by
``compressing'' the size of the matrix for the pivoted QR
algorithm. In our previous work \cite{LuYing:15}, we apply a random
Fourier projection on $\rho$, which amounts to taking a sample of
random linear combinations of rows of the $\rho$ matrix:
\begin{equation}\label{eq:originalfourierproj}
  \varrho_{\xi}(x) = \sum_{\gamma} e^{-i 2\pi \gamma \xi (N^2K^2)^{-1}} \eta_{\gamma} \rho_{\gamma}(x)
\end{equation}
for random samples of $\xi \in \{1, \ldots, N^2K^2\}$ of total number
proportional to $N$ and $\eta_{\gamma}$ being random unit complex
number for each $\gamma$. Here, we have used $\gamma \in \{1, \ldots,
N^2K^2\}$ to denote the row index $(nkml)$ for simplicity of
notations. Fixing $\xi$, $\varrho_{\xi}$ is a random linear
combination of the rows $\rho_{\gamma}$. Thus, important columns of
$\rho$ remains important columns of $\varrho$ with high probability.

The total size of the sample $\xi$ to guarantee successful
identification of important columns is $\Or(N)$ and is much smaller
than the original number of rows of $\rho$, and hence the cost of the
pivoted QR algorithm is much reduced. The computational cost is
dominated by using FFT to compute \eqref{eq:originalfourierproj}, and
is hence $\Or\bigl(N_{\grid} (NK)^2 \log (NK)\bigr)$.

While the efficiency is much improved by using the algorithm in
\cite{LuYing:15}, it is still expensive in the current context due to
the large number of rows of the matrix $\rho$. Our key observation to
further reduce the computational cost is that $\rho$ given by
$\rho_{nkml}(x) = \wb{u}_{n,k}(x) u_{m,l}(x)$ has the tensor product
structure, in the sense that the row index $(nkml)$ is naturally
separated into two groups $(nk)$ and $(ml)$. We may explore this
structure by taking the random Fourier projection on $u$ instead,
which has a much smaller dimension than $\rho$, and hence makes the
row projection more efficient. More concretely, we consider instead a
projection of the type (cf. \eqref{eq:originalfourierproj})
\begin{equation}\label{eq:newfourierproj}
  M_{ij}(x) = \sum_{\alpha} e^{i 2\pi \alpha \xi_i  (NK)^{-1}} \wb{\eta}_{\alpha} \wb{u}_{\alpha}(x) 
  \sum_{\beta} e^{-i 2\pi \beta \xi_j  (NK)^{-1}} \eta_{\beta} u_{\beta}(x) 
\end{equation}
for $\{\xi_i\} \subset  \{1, \ldots, NK\}$ and $\alpha, \beta \in \{1,
\ldots, NK\}$ being the row index of $u$, viewed as a $(NK) \times
N_{\grid}$ matrix. Note that we can rewrite \eqref{eq:newfourierproj}
as
\begin{equation}
  \begin{aligned}
  M_{ij}(x) = \sum_{\alpha, \beta} e^{i 2\pi (\alpha \xi_i - \beta \xi_j) (NK)^{-1}} 
  \wb{\eta}_{\alpha} \eta_{\beta} \wb{u}_{\alpha}(x)  u_{\beta}(x) 
   = \sum_{\alpha, \beta} e^{i 2\pi (\alpha \xi_i - \beta \xi_j) (NK)^{-1}} 
  \wb{\eta}_{\alpha} \eta_{\beta} \rho_{\alpha\beta}(x).
  \end{aligned}
\end{equation}
and hence this still corresponds to a random linear combination of the
rows of $\rho$, and hence could be still used to select important
columns of $\rho$ to represent the column space. The details of the
algorithm are described in Algorithm~\ref{alg:df}.

\begin{algorithm2e}[H]
  \SetKwInOut{Input}{Input}
  \SetKwInOut{Output}{Output}

  \Input{Orbitals $u_{nk}(x)$} 
  \Output{Selected grids $x_{\mu}$ and
    auxiliary basis functions $P_{\mu}(x)$, such that $\wb{u}_{n,k}(x)
    u_{m,l}(x) \approx \sum_{\mu} \wb{u}_{n,k}(x_{\mu})
    u_{m,l}(x_{\mu}) P_{\mu}(x)$}

  Reshape $u_{nk}(x)$ into an $(NK) \times N_{\grid}$ matrix
  $U_{\alpha}(x)$ where $\alpha = 1, \ldots, NK$ is the row index;

  Compute the discrete Fourier transform of $U$ left multiplied by
  a random diagonal matrix:
  \begin{equation*}
    \wh{U}_{\xi}(x) = \sum_{\alpha} e^{-i 2\pi \alpha \xi  (NK)^{-1}} \eta_{\alpha} U_{\alpha}(x), 
  \end{equation*}
  where $\eta_{\alpha}$ is a random unit complex number for each $\alpha$; 

  Choose a submatrix $\mc{U}$ of $\wh{U}$ by randomly choosing
  $r = c \sqrt{N}$ rows.

  Construct a $r^2 \times N_{\grid}$ matrix $M$  
  \begin{equation*}
    M_{ij}(x) = \wb{\mc{U}}_{i}(x) \mc{U}_{j}(x), \qquad i,j = 1, \ldots, r,
  \end{equation*}
  where we view $(ij)$ as the row index of $M$. 

  Apply column selection Algorithm~\ref{alg:cs} on the $r^2 \times
  N_{\grid}$ matrix $M$ to find selected $N_\col$ columns
  $\{x_{\mu}\}$ and auxiliary basis functions $P_{\mu}(x)$.

  \caption{Density fitting via random Fourier projection and column
    selection} \label{alg:df}
\end{algorithm2e}

The computationally costly steps in Algorithm~\ref{alg:df} are:
\begin{enumerate}
\item Step 2 requires FFT of vectors of length $(NK)$ for $N_{\grid}$
  times, and hence has computational cost $\Or\bigl(N_{\grid} NK \log
  (NK)\bigr)$;
\item Step 4 scales as $\Or(r^2 N_{\grid}) = \Or(N N_{\grid})$ by the choice of $r = c \sqrt{N}$; 
\item Step 5 requires pivoted QR algorithm applied on a $N \times
  N_{\grid}$ matrix, and hence has cost $\Or(N^2 N_{\grid})$, which dominates the other steps in Algorithm~\ref{alg:cs}. 
\end{enumerate}
Therefore, the total computational cost of the algorithm scales as
$\Or\bigl(N_{\grid} N^2 + N_{\grid} NK \log (NK)\bigr)$. In practice,
the prefactor of the FFT is much smaller than the one of the pivoted
QR factorization. As a result, the $\Or\bigl(N_{\grid} N^2\bigr)$ part
dominates the actual running time.

\section{Numerical examples}\label{sec:numerics}

In this section, we consider a few model examples to demonstrate the
effectiveness of the algorithms proposed above.  In order to measure
the error of the density fitting, two metrics can be used. The error
in $L^2$ metric measures
\begin{equation}
  \norm{\rho_{nkml} - \wt{\rho}_{nkml}}_2, 
\end{equation}
for $\wt{\rho}_{nkml} = \sum_{\mu} \wb{C}_{nk}^{\mu} C_{ml}^{\mu} P_{\mu}$, 
while the error in Coulomb metric is given by 
\begin{equation}
  \norm{\rho_{nkml} - \wt{\rho}_{nkml}}_C = \Biggl( \iint_{\Gamma \times \Gamma}
  \bigl(\rho_{nkml} - \wt{\rho}_{nkml}\bigr)(x) G(x - y) 
    \bigl(\rho_{nkml} - \wt{\rho}_{nkml}\bigr)(y)  \ud x \ud y \Biggr)^{1/2}, 
\end{equation}
where $G$ is the periodic Coulomb kernel, solving 
\begin{equation}
  -\Delta G( \cdot - y) = 4\pi ( \delta_{y} - 1)
\end{equation}
with periodic boundary condition on $\Gamma = [0, 1]^d$ and
$\int_{\Gamma} G = 0$.

The Coulomb metric is a good measure of the error since density
fitting is often used in approximating electron repulsion integral
tensor, such as
\begin{equation}
  \begin{aligned}
    E_{nkml} & = \iint_{\Gamma \times \Gamma} \wb{\psi}_{n, k}(x)
    \psi_{m, l}(x)
    G(x - y) \wb{\psi}_{m, l}(y) \psi_{n, k}(y) \ud x \ud y \\
    & = \iint_{\Gamma \times \Gamma} \wb{u}_{n, k}(x) u_{m, l}(x) G(x
    - y) \wb{u}_{m, l}(y) u_{n, k}(y) e^{- i (k - l) \cdot (x - y)}
    \ud x\ud y.
  \end{aligned}
\end{equation}
Using the density fitting \eqref{eq:dfapprox}, we get the
approximation
\begin{equation}
  \begin{aligned}
    \wt{E}_{nkml} & = \sum_{\mu\nu} \iint_{\Gamma \times \Gamma}
    \wb{C}_{n,
      k}^{\mu} C_{m, l}^{\mu} P_{\mu}(x)  G(x-y) \wb{C}_{m, l}^{\nu} C_{n, k}^{\nu} P_{\nu}(y) e^{- i (k - l) \cdot (x - y)}  \ud x \ud y \\
    & = \sum_{\mu\nu} \wb{C}_{m, k}^{\mu} C_{n, l}^{\mu} \wb{C}_{m,
      l}^{\nu} C_{n, k}^{\nu} \iint_{\Gamma\times \Gamma} P_{\mu}(x)
    G(x-y) P_{\nu}(y) e^{-i (k-l)\cdot(x-y)} \ud x \ud y.
  \end{aligned}
\end{equation}
A simple calculation gives 
\begin{equation}
  \begin{aligned}
    \abs{E_{nkml} - \wt{E}_{nkml}} & = \Biggl\lvert \iint \Bigl(
    \rho_{nkml}(x) \wb{\rho}_{nkml}(y) - \wt{\rho}_{nkml}(x)
    \wb{\wt{\rho}}_{nkml}(y) \Bigr) G(x - y) e^{-i(k-l) \cdot (x - y)}
    \ud x \ud y\,  \Biggr\rvert \\
    & \leq \Biggl\lvert \iint \bigl( \rho_{nkml}(x) -
    \wt{\rho}_{nkml}(x) \bigr) \wb{\rho}_{nkml}(y) G(x - y) e^{-i(k-l)
      \cdot (x - y)} \ud x \ud y\,  \Biggr\rvert \\
    & \qquad + \Biggl\lvert \iint \wt{\rho}_{nkml}(x) \bigl(
    \wb{\rho}_{nkml}(y) - \wb{\wt{\rho}}_{nkml}(y) \bigr) G(x - y)
    e^{-i(k-l) \cdot (x - y)}  \ud x \ud y\,  \Biggr\rvert \\
    & \leq C \norm{\rho_{nkml} - \wt{\rho}_{nkml}}_C \bigl(
    \norm{\rho_{nkml}}_C + \norm{\wt{\rho}_{nkml}}_C \bigr),
  \end{aligned}
\end{equation}
where the last step uses the Fourier representation and the 
Cauchy-Schwartz inequality.

\subsection{2D examples}

We consider two examples in two dimensions. In both examples, the
periodic model potential $V(x)$ is given by
\[
V(x) = \sum_{n_1\in \ZZ}\sum_{n_2\in \ZZ} V_c(x-n_1 e_1-n_2 e_2)
\]
where $V_c$ is a localized potential centered at the origin and $e_1$
and $e_2$ are the unit Cartesian basis vectors. This is a simple
square lattice and the Brillouin zone is also a square. The potential
$V_c(\cdot)$ is chosen to be centered at the origin and spherically
symmetric. The unit cell $[0,1)^2$ is discretized with $N_{\grid}=48$
  points in each dimension. In each example, we consider up to $N=41$
  bands and the Brillouin zone $[-\pi,\pi)^2$ is sampled with up to
    $K=16$ points per dimension. The band structure and the
    eigenfunctions are computed with spectral discretization for high
    order accuracy and the resulting discrete eigenvalue problem for
    each $k$ point is solved using LOBPCG \cite{Knyazev:01} with the
    simple inverse Laplacian preconditioning. The prescribed accuracy
    for the column selection algorithm is set to be of order
    $10^{-5}$.

\begin{figure}[ht!]
  \begin{tabular}{cc}
    \includegraphics[height=1.5in]{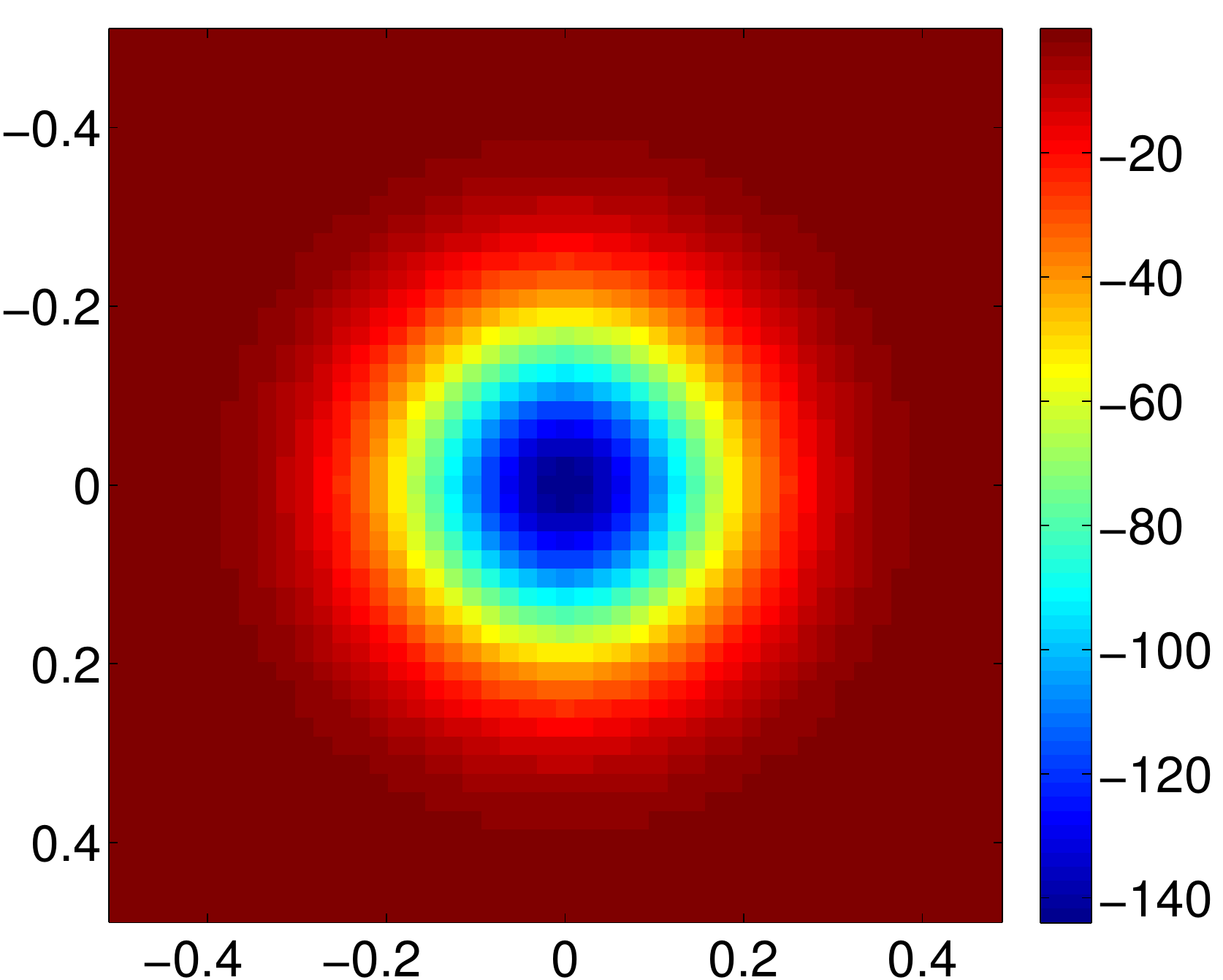} & \includegraphics[height=1.5in]{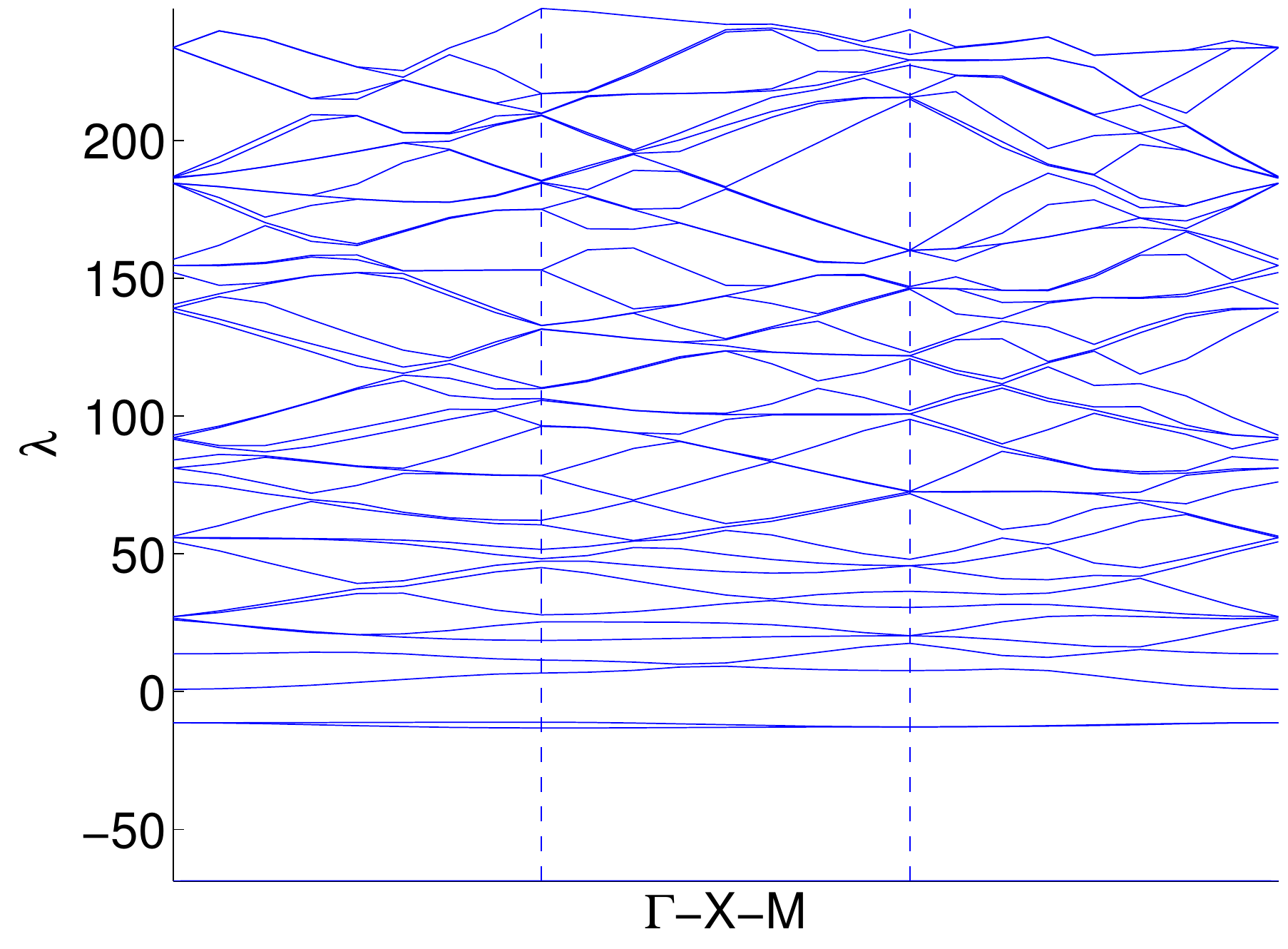} \\
    (a) & (b) \\
    \includegraphics[height=1.5in]{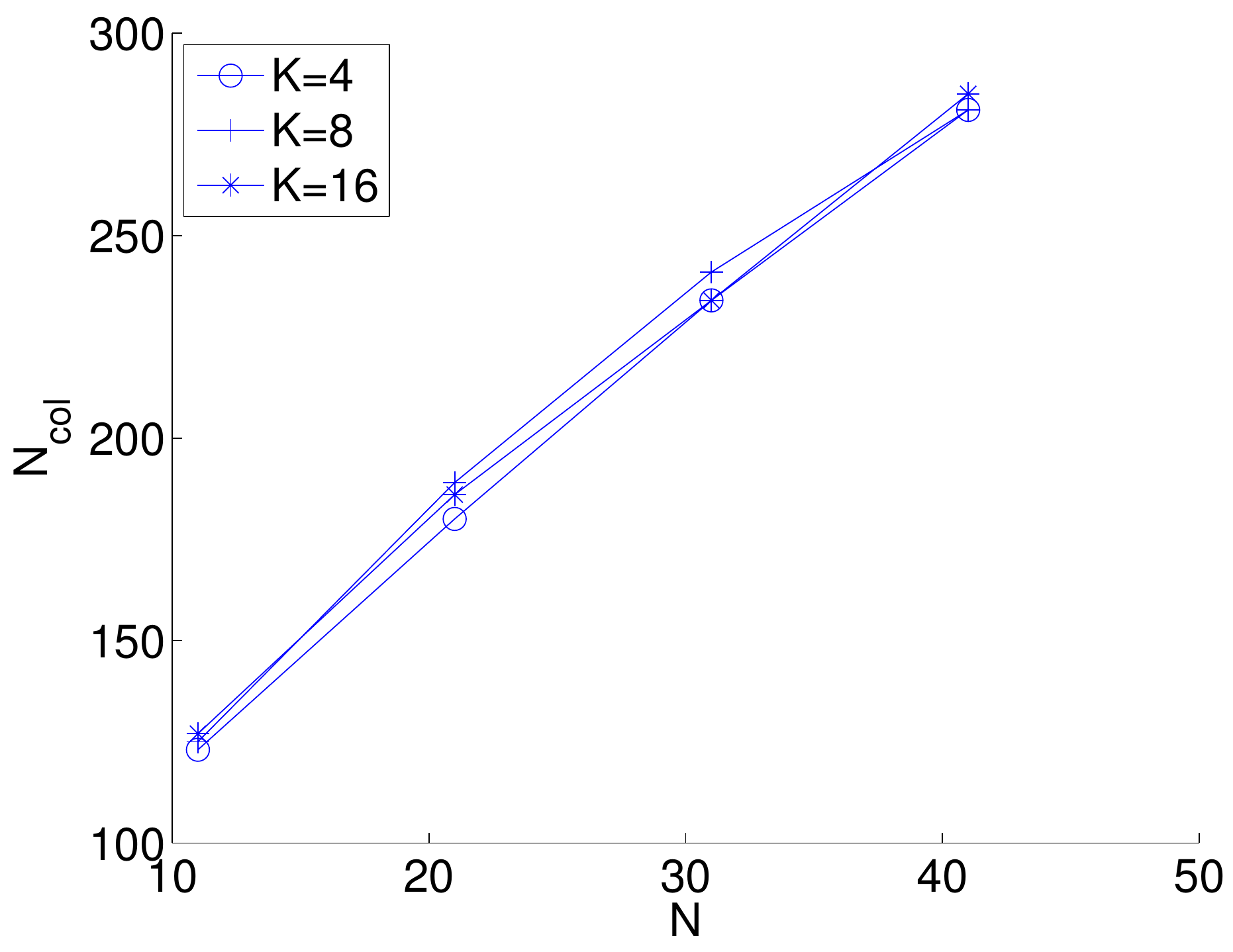} & \includegraphics[height=1.5in]{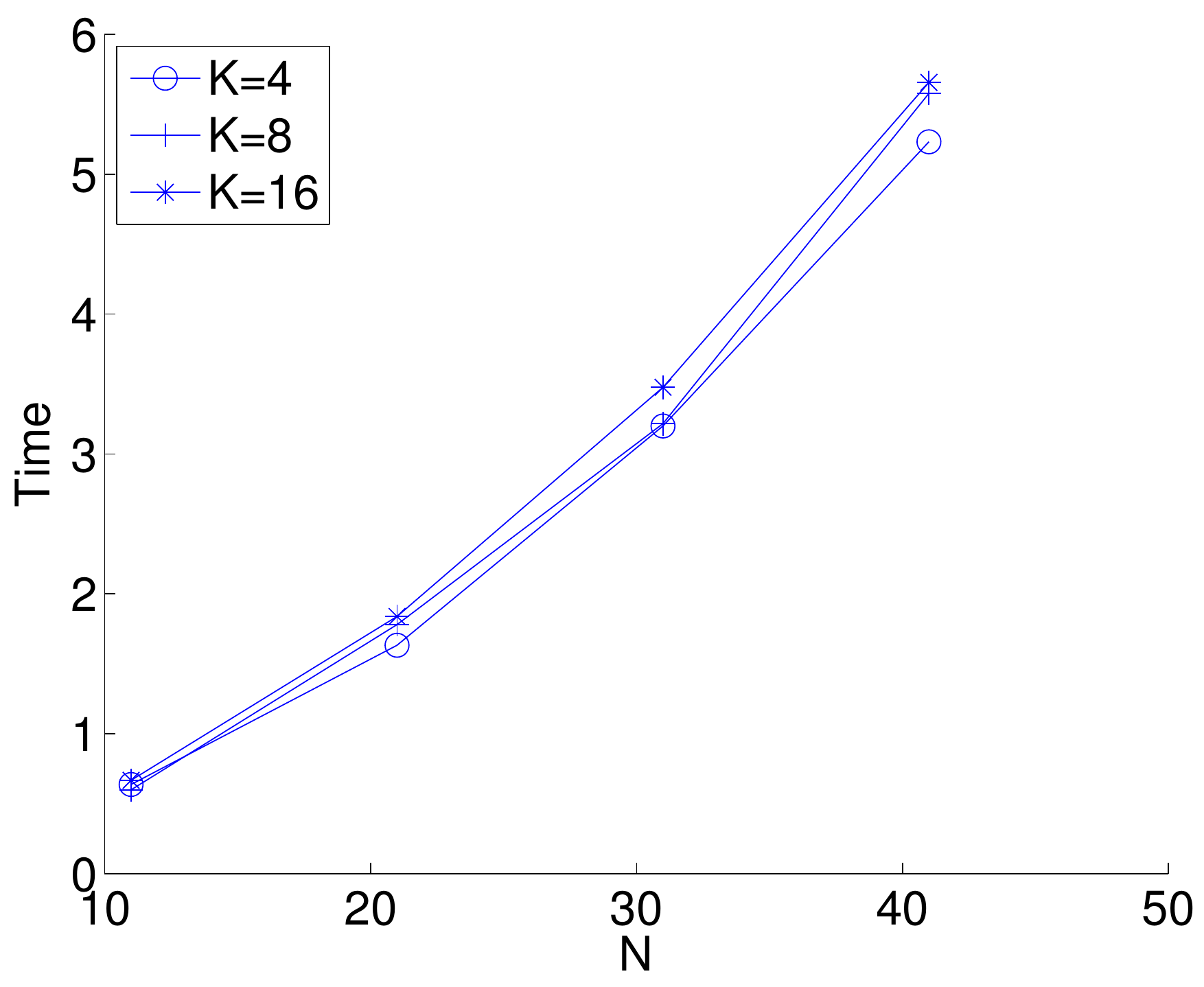} \\
    (c) & (d) \\
    \includegraphics[height=1.5in]{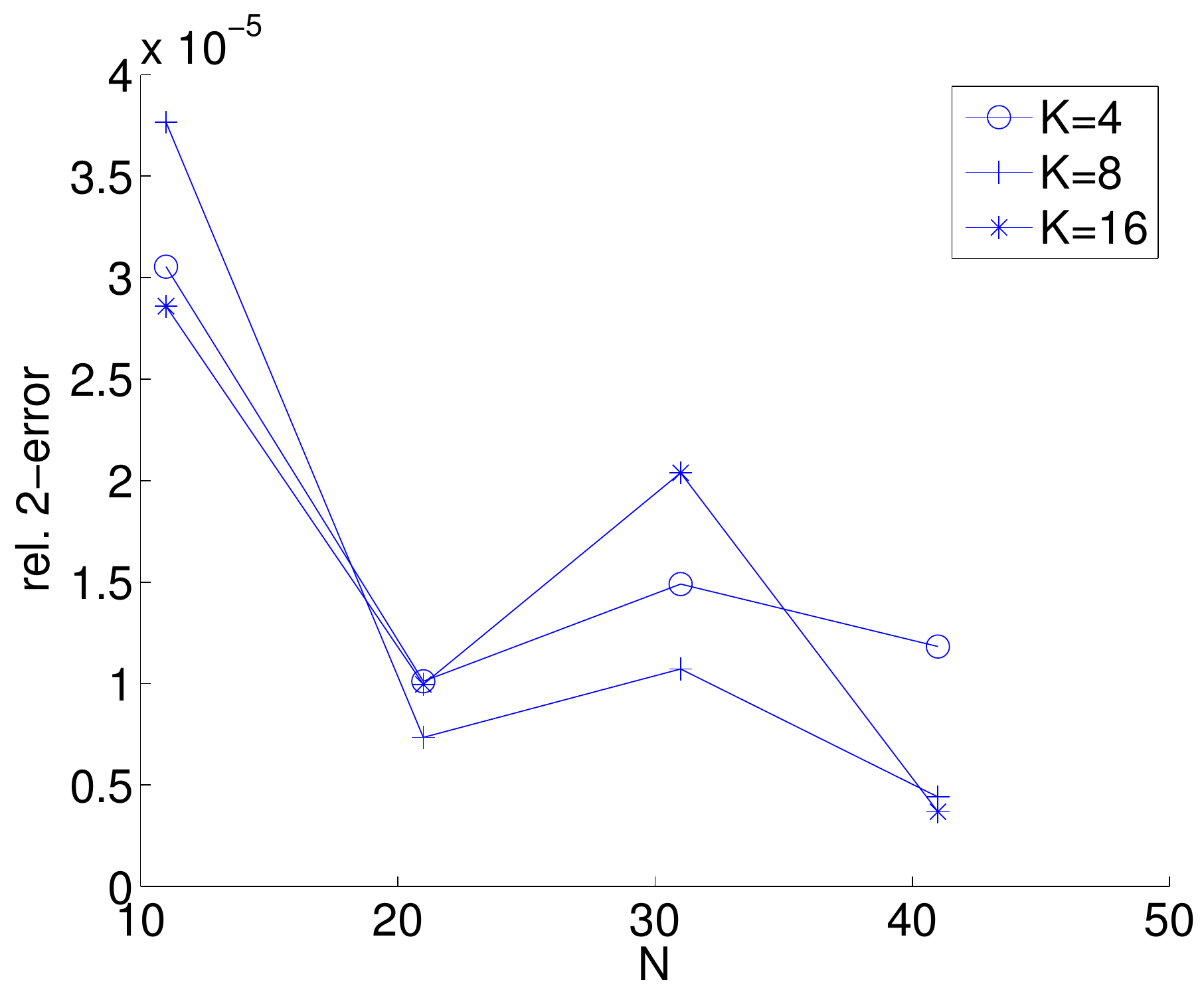} & \includegraphics[height=1.5in]{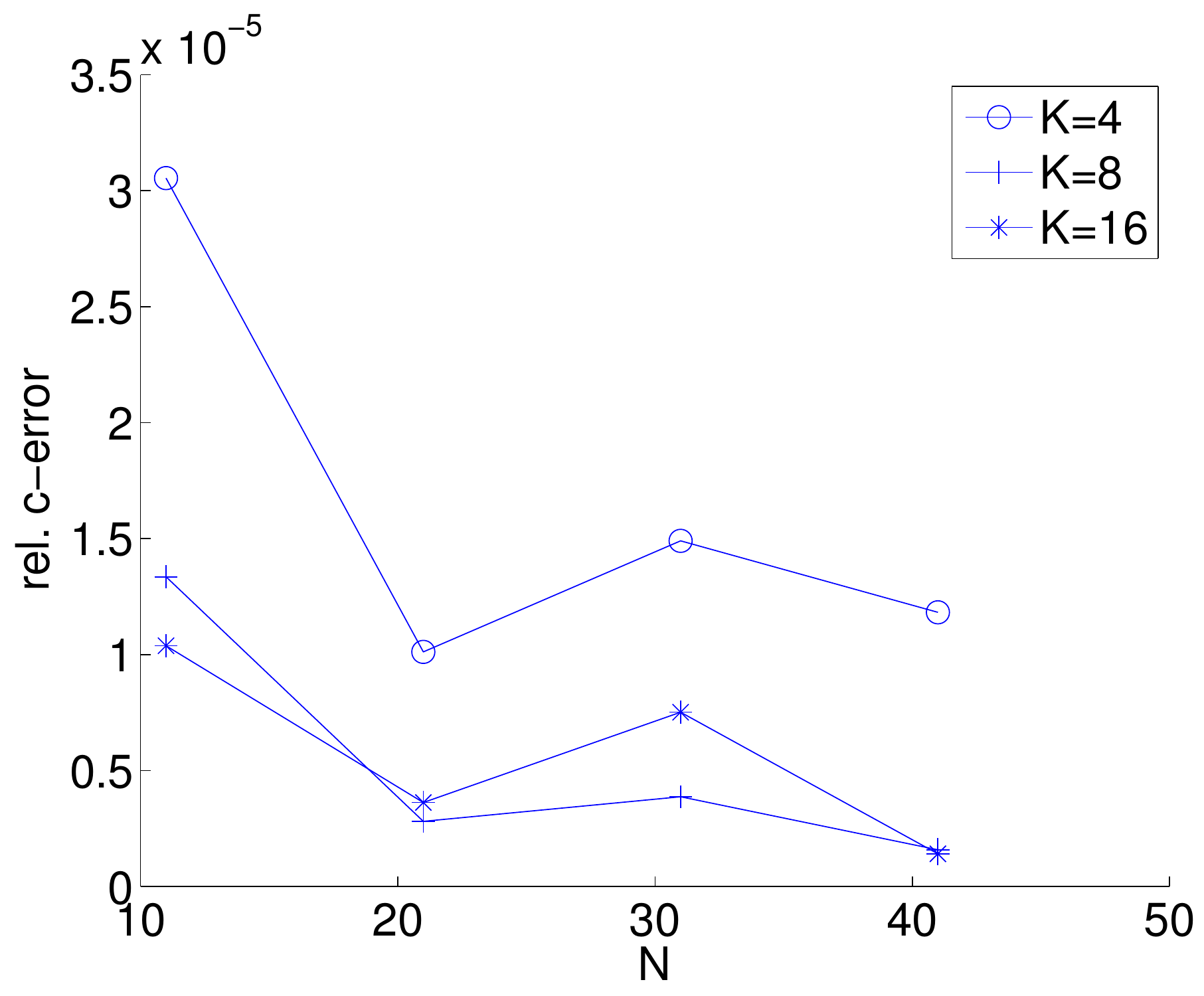}\\
    (e) & (f)
  \end{tabular}
  \caption{Results of the first 2D example.
    (a) the potential $V_c$ in unit cell.
    (b) the band structure plotted along the $\Gamma-X-M$ path.
    (c) $N_{\col}$ as a function $N$ (the number of bands) for different values of $K$ (the number of $k$ points).
    (d) The time used by the column selection algorithm in seconds.
    (e) The relative error measured in the $L^2$ metric.
    (f) The relative error measured in the Coulomb metric.
  }
  \label{fig:t2d4}
\end{figure}

In the first example, the potential $V_c(x)$ is given by a Gaussian
profile 
\begin{equation}\label{eq:t2d4}
  V_c(x) = -144 \exp\left(-\frac{\|x\|^2}{2\sigma^2}\right)
\end{equation}
with $\sigma=0.1333$.  Figure \ref{fig:t2d4} summarizes the results for
this first example.
\begin{itemize}
\item Figure \ref{fig:t2d4}(a) shows the potential profile in the unit
  cell $[-1/2,1/2)^2$. 
\item Figure \ref{fig:t2d4}(b) gives the band structure along the
  standard $\Gamma-X-M$ path in the Brillouin zone for simple square
  lattice. 
\item Figure \ref{fig:t2d4}(c) plots the dependence of $N_{\col}$ (the
  number of auxiliary basis) as a function of $N$ (the number of
  bands) for different choices of $K$ (the number of $k$-point samples
  in each dimension) at the targeted accuracy level $10^{-5}$. From
  this plot, one can clearly see that $N_{\col}$ grows roughly
  linearly with respect to $N$ and that it depends only mildly on the
  number of $k$ point samples.
\item The time for performing the column selection in seconds is given
  in Figure \ref{fig:t2d4}(d) as a function of $N$ for different
  values of $K$. This plot shows that the complexity depends
  quadratically on $N$ but only grows mildly with $K$.
\item Finally, Figures \ref{fig:t2d4}(e) and (f) give the relative
  errors in the $L^2$ metric and the Coulomb metric,
  respectively. These two plots show that the estimated relative
  errors are bounded by a small constant times the prescribed accuracy
  level.
\end{itemize}

\begin{figure}[ht!]
  \begin{tabular}{cc}
    \includegraphics[height=1.5in]{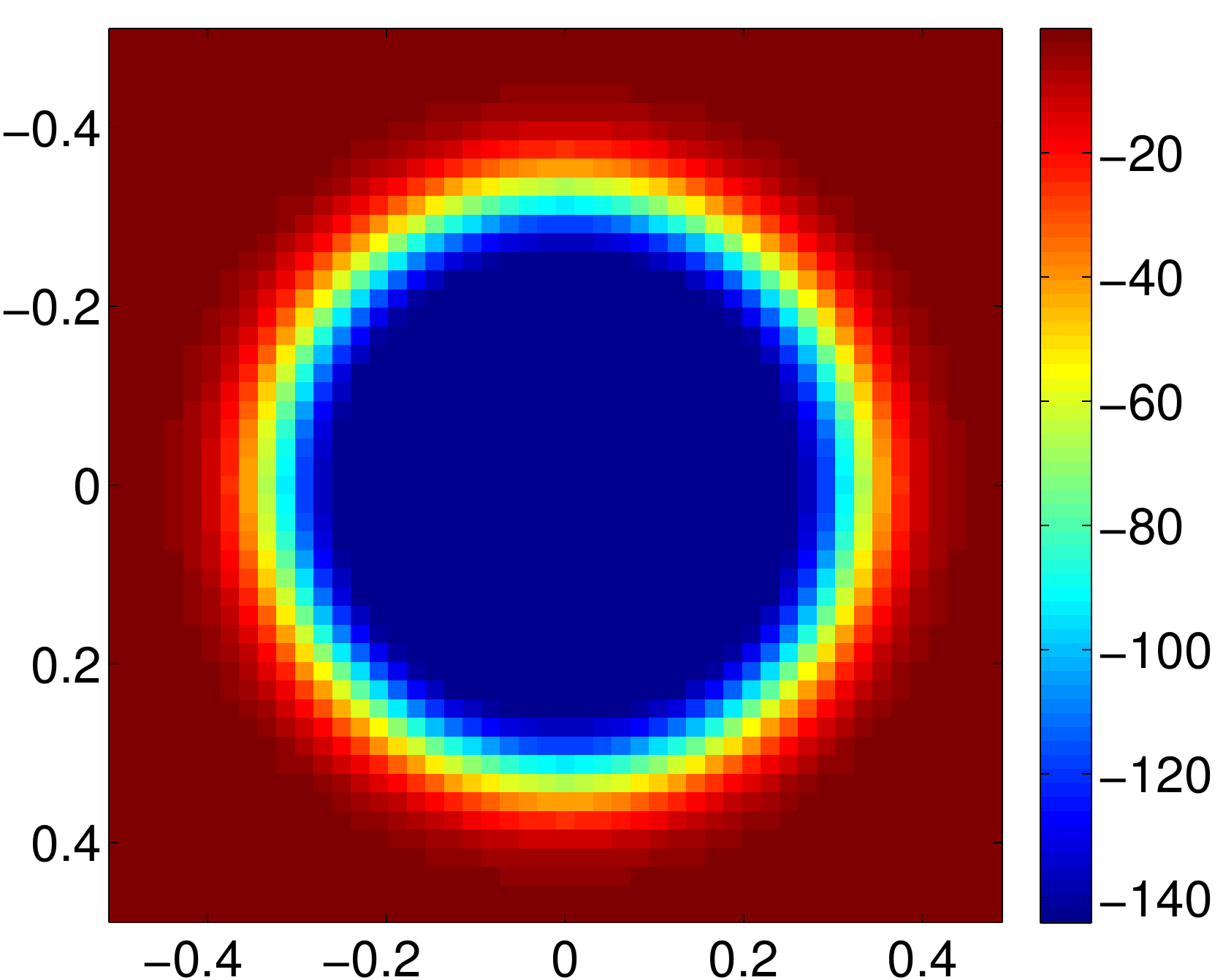} & \includegraphics[height=1.5in]{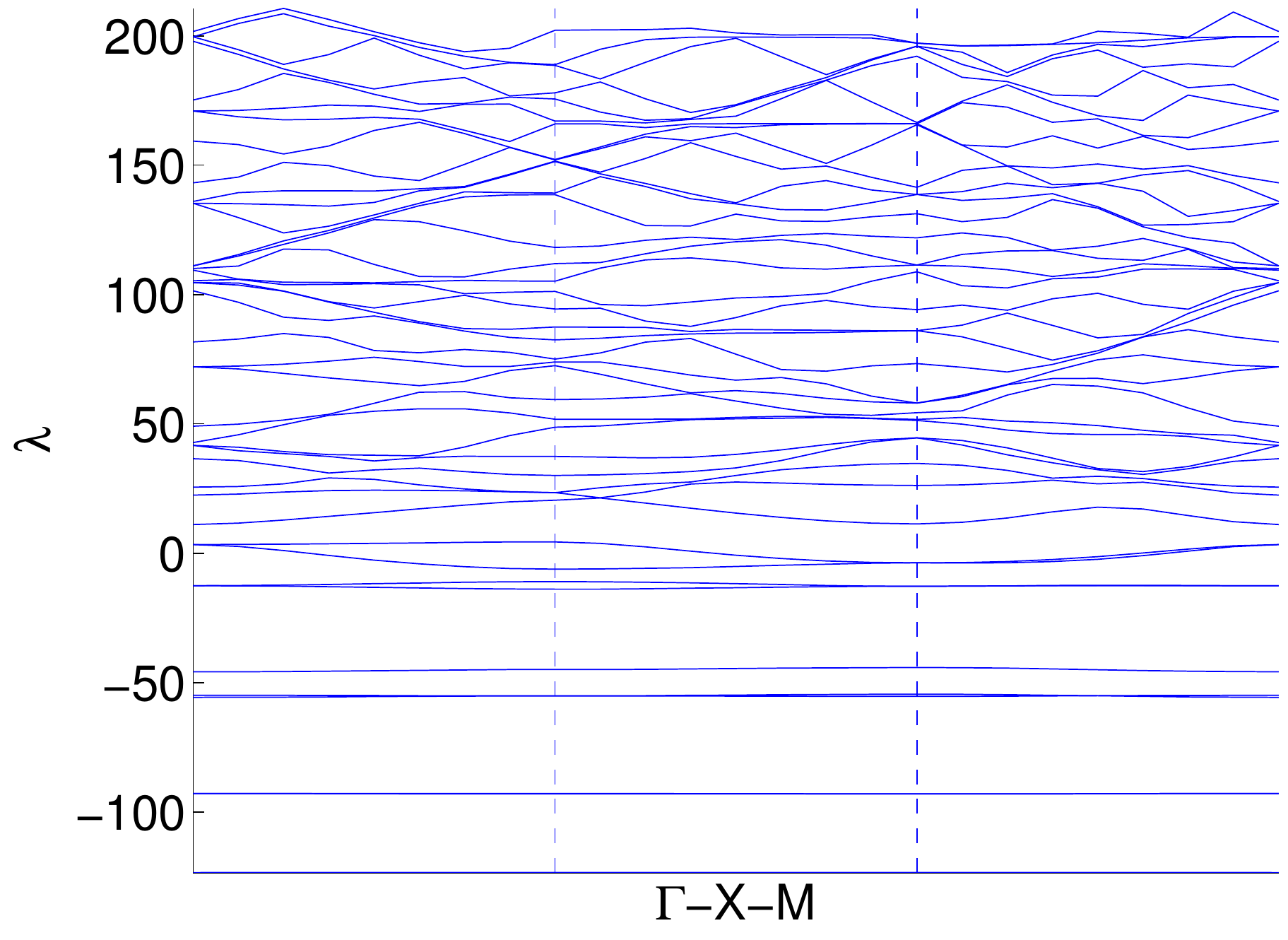} \\
    (a) & (b) \\
    \includegraphics[height=1.5in]{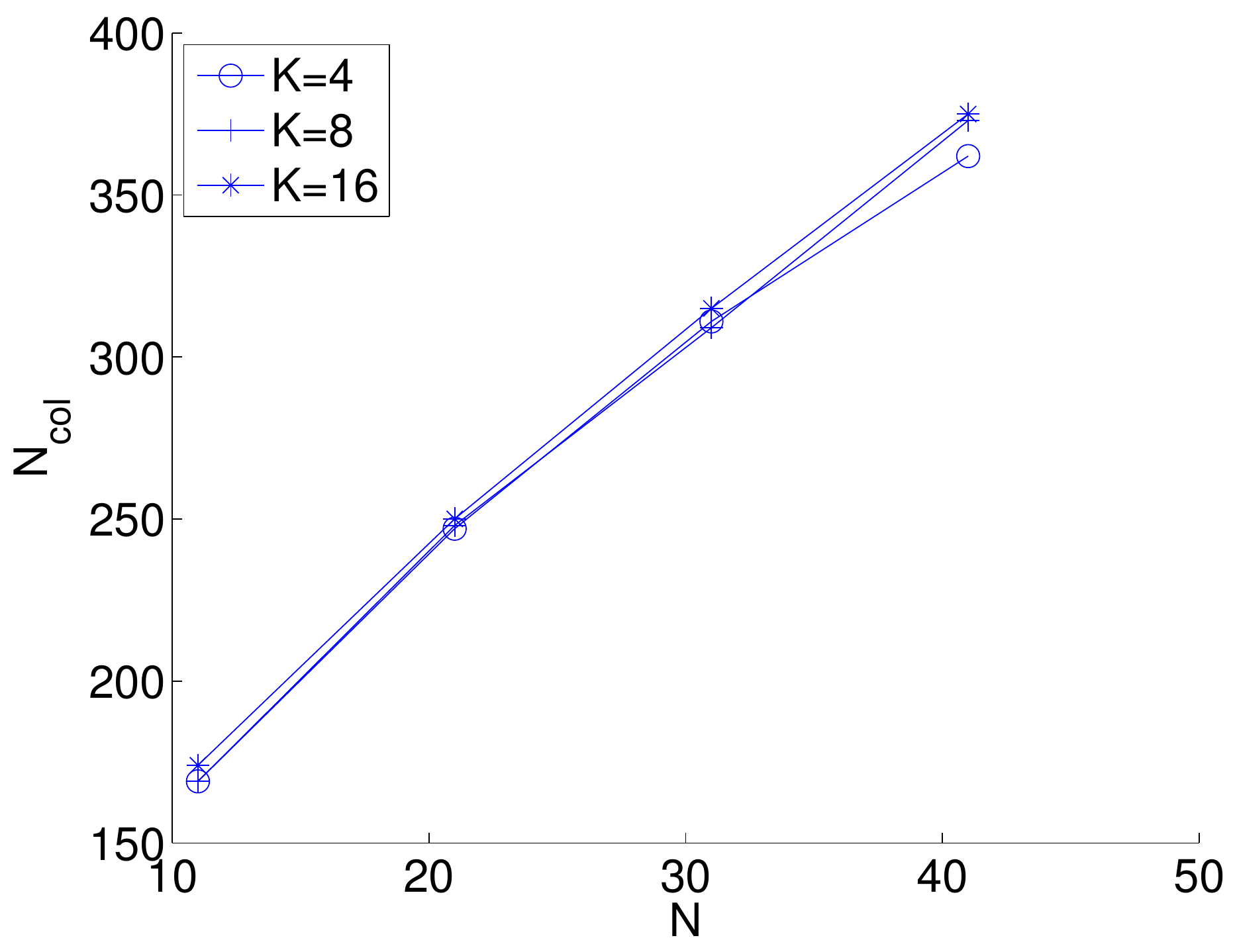} & \includegraphics[height=1.5in]{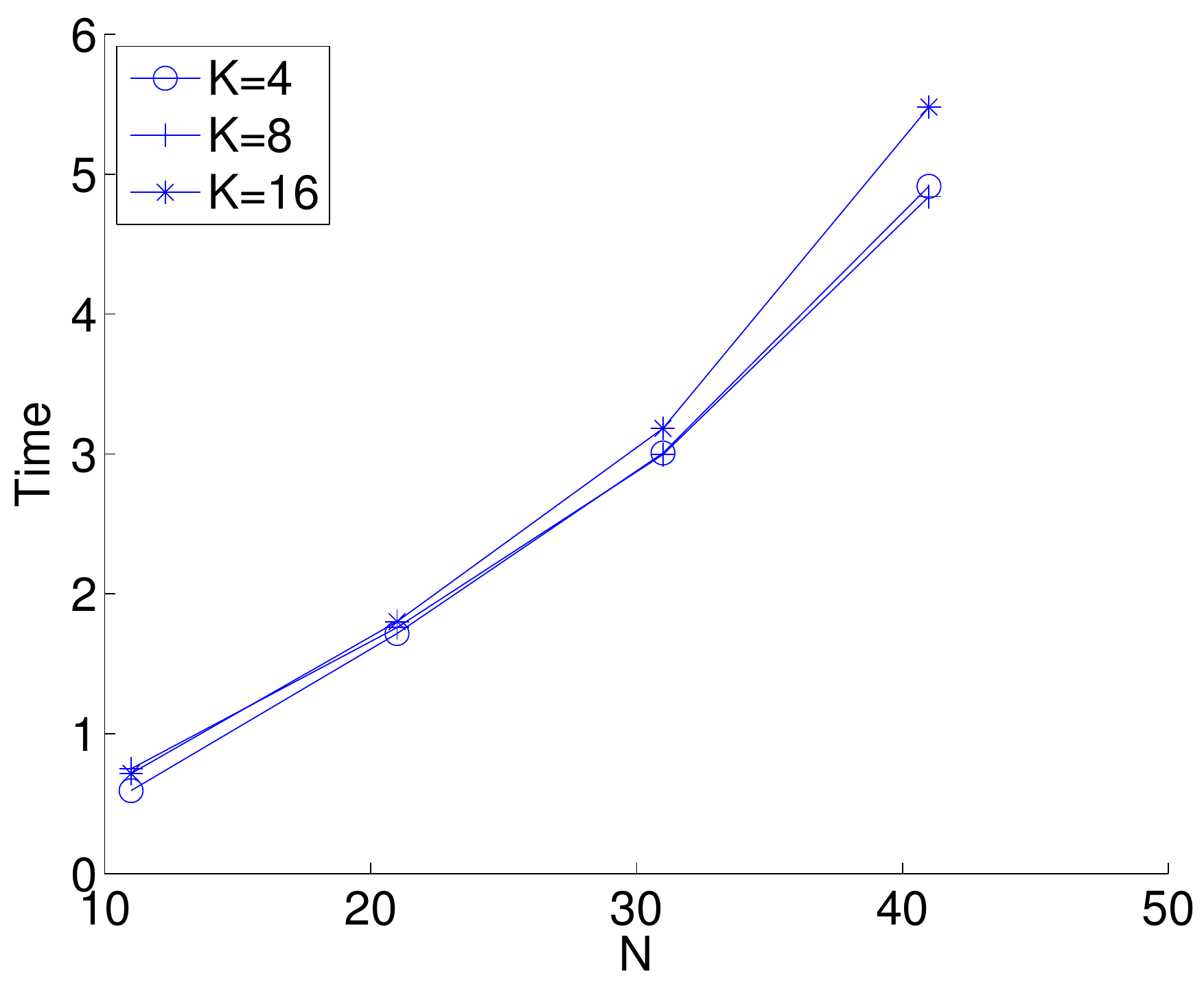} \\
    (c) & (d) \\
    \includegraphics[height=1.5in]{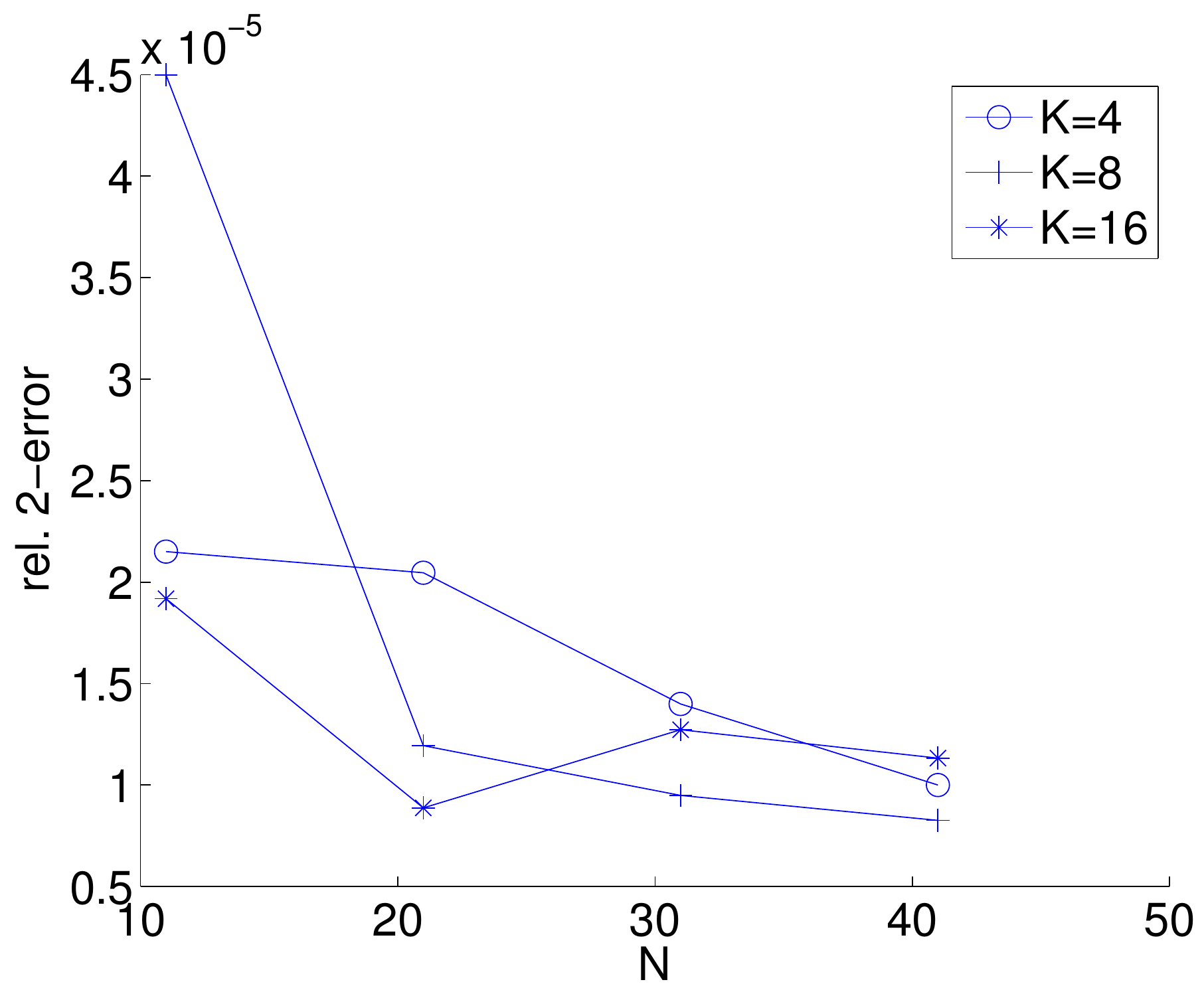} & \includegraphics[height=1.5in]{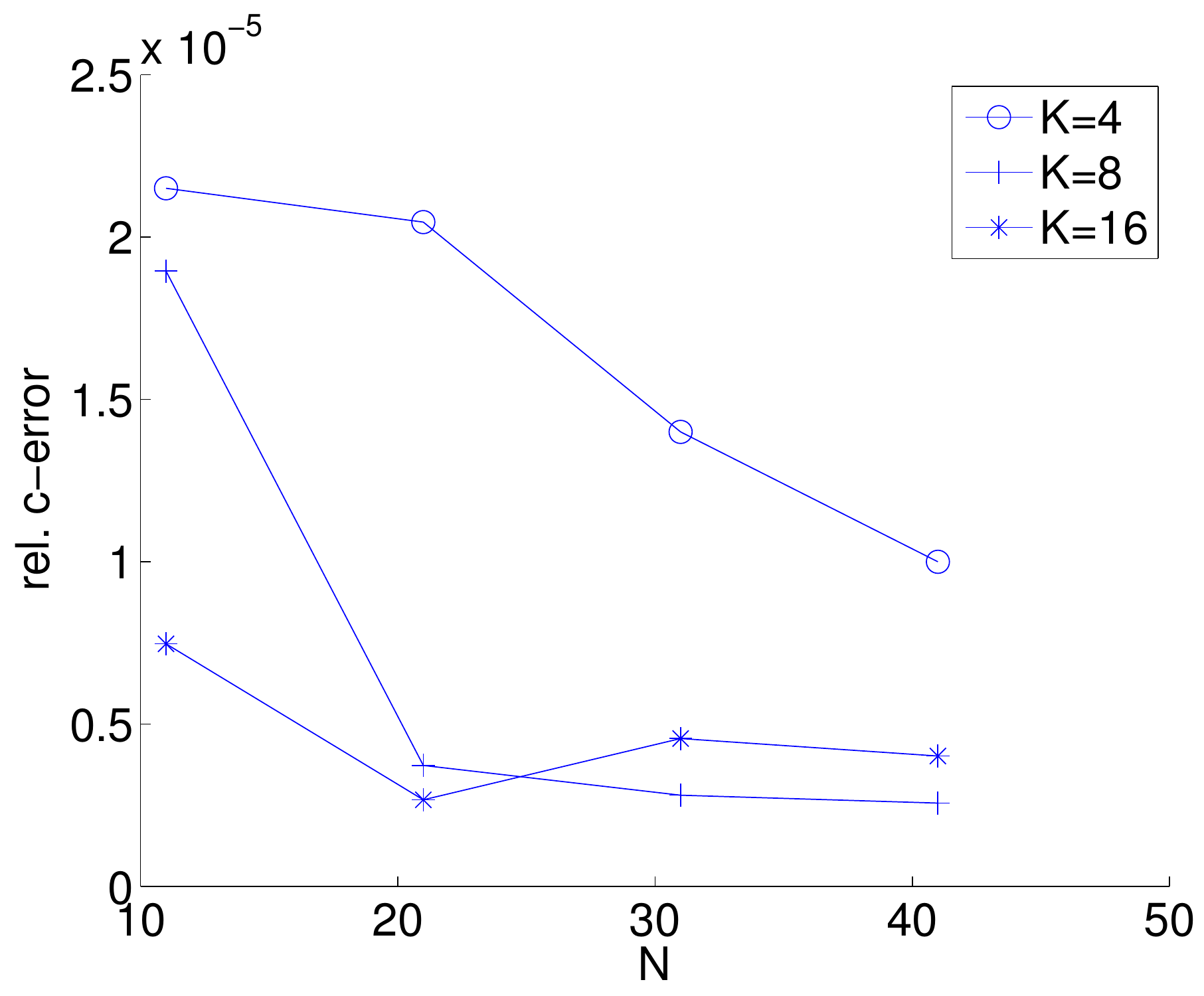}\\
    (e) & (f)
  \end{tabular}
  \caption{Results of the second 2D example.
    (a) the potential $V_c$ in unit cell.
    (b) the band structure plotted along the $\Gamma-X-M$ path.
    (c) $N_{\col}$ as a function $N$ (the number of bands) for different values of $K$ (the number of $k$ points).
    (d) The time used by the column selection algorithm in seconds.
    (e) The relative error measured in the $L^2$ metric.
    (f) The relative error measured in the Coulomb metric.
  }
  \label{fig:t2d5}
\end{figure}

In the second example, the potential $V_c(x)$ is chosen to be 
\begin{equation}\label{eq:t2d5}
  V_c(x) = -144 \exp\left( -\frac{\max(\|x\|-1/4,0)^2}{2\sigma^2} \right)
\end{equation}
with $\sigma = 0.0667$. The results for this example are summarized in
Figure \ref{fig:t2d5}. The results show similar characteristics as the
previous example.

\subsection{3D examples}

For the 3D case, we consider two similar examples. In each example,
the periodic model potential $V(x)$ is given by
\[
V(x) = \sum_{n_1\in \ZZ}\sum_{n_2\in \ZZ}\sum_{n_3\in \ZZ} V_c(x-n_1 e_1-n_2 e_2-n_3 e_3)
\]
where $V_c$ is a localized potential centered at the origin and $e_1$,
$e_2$, and $e_3$ are the unit Cartesian basis vectors. This is a
simple cubic lattice and the Brillouin zone is also a cube. The
potential $V_c(\cdot)$ is again chosen to be centered at the origin
and spherically symmetric. The unit cell $[0,1)^3$ is discretized with
  $N_\grid=24$ points in each dimension. In each example, we consider
  up to $N=41$ bands and the Brillouin zone $[-\pi,\pi)^3$ is sampled
    with up to $K=12$ points per dimension. Similar to the 2D
    examples, the band structure and the eigenfunctions are computed
    with pseudo-spectral discretization for high order accuracy and
    the resulting discrete eigenvalue problem for each $k$ point is
    solved using LOBPCG with simple inverse Laplacian
    preconditioning. The prescribed accuracy for the column selection
    algorithm is set to be of order $10^{-5}$.

\begin{figure}[ht!]
  \begin{tabular}{cc}
    \includegraphics[height=1.5in]{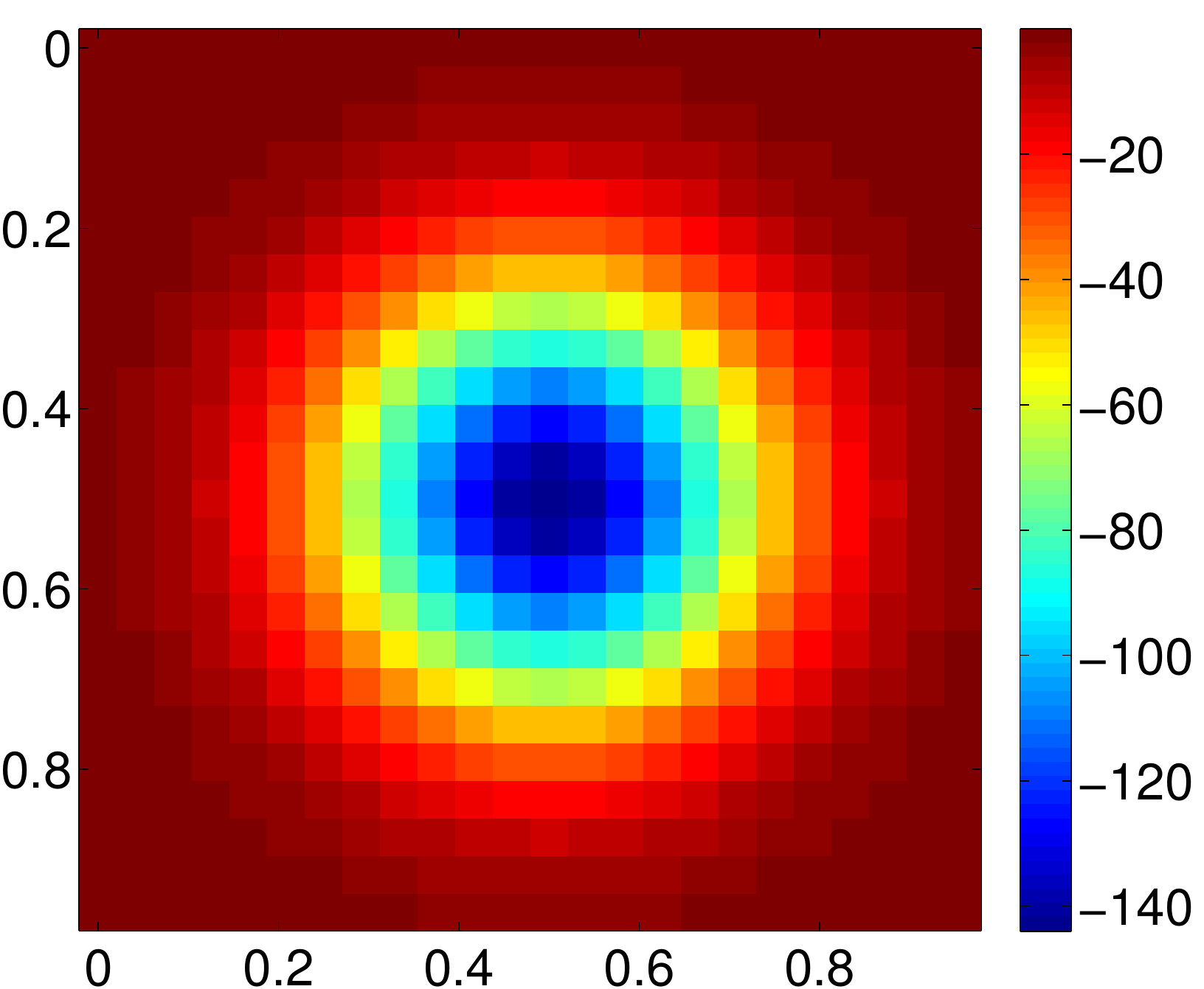} & \includegraphics[height=1.5in]{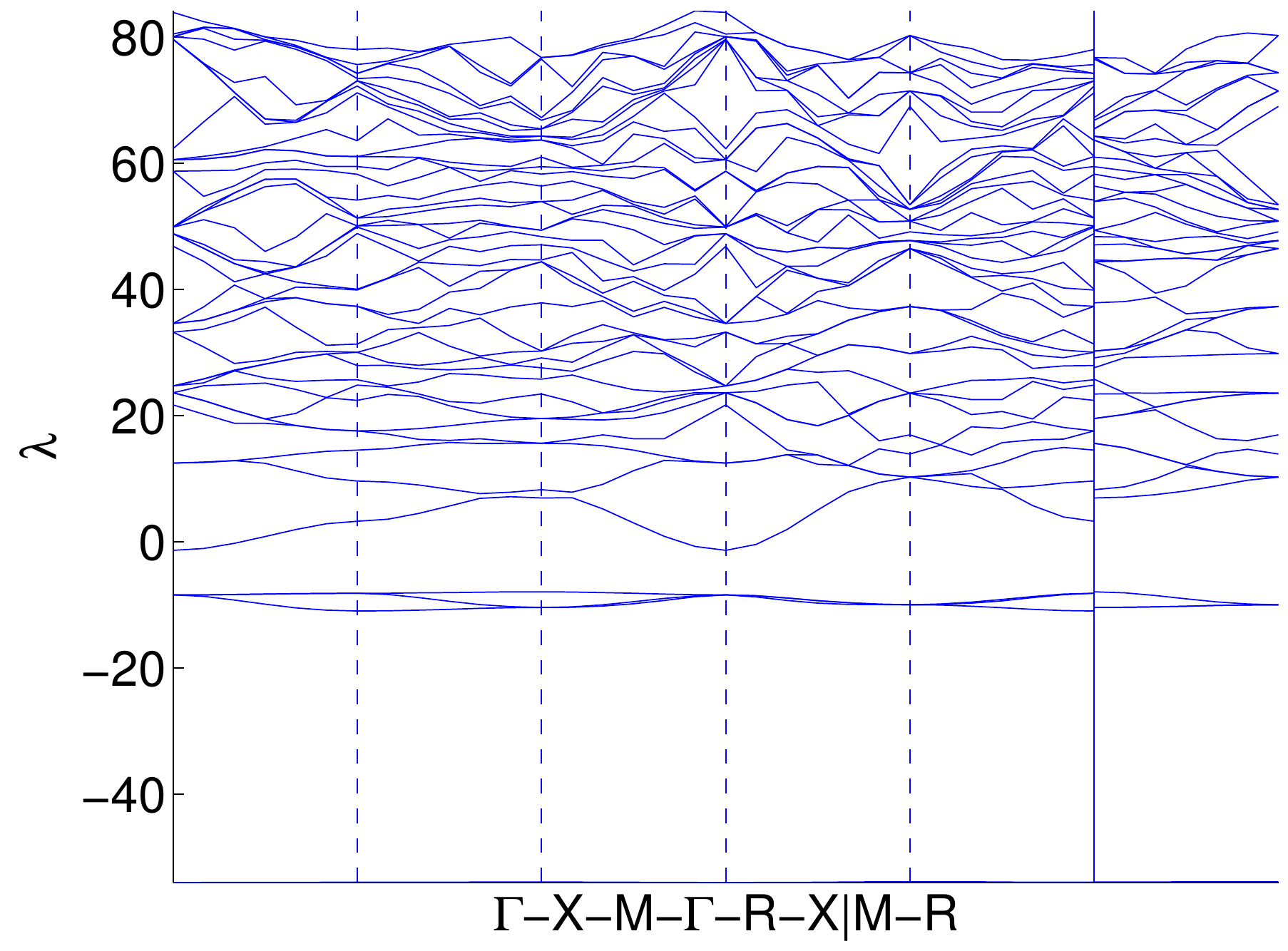} \\
    (a) & (b) \\
    \includegraphics[height=1.5in]{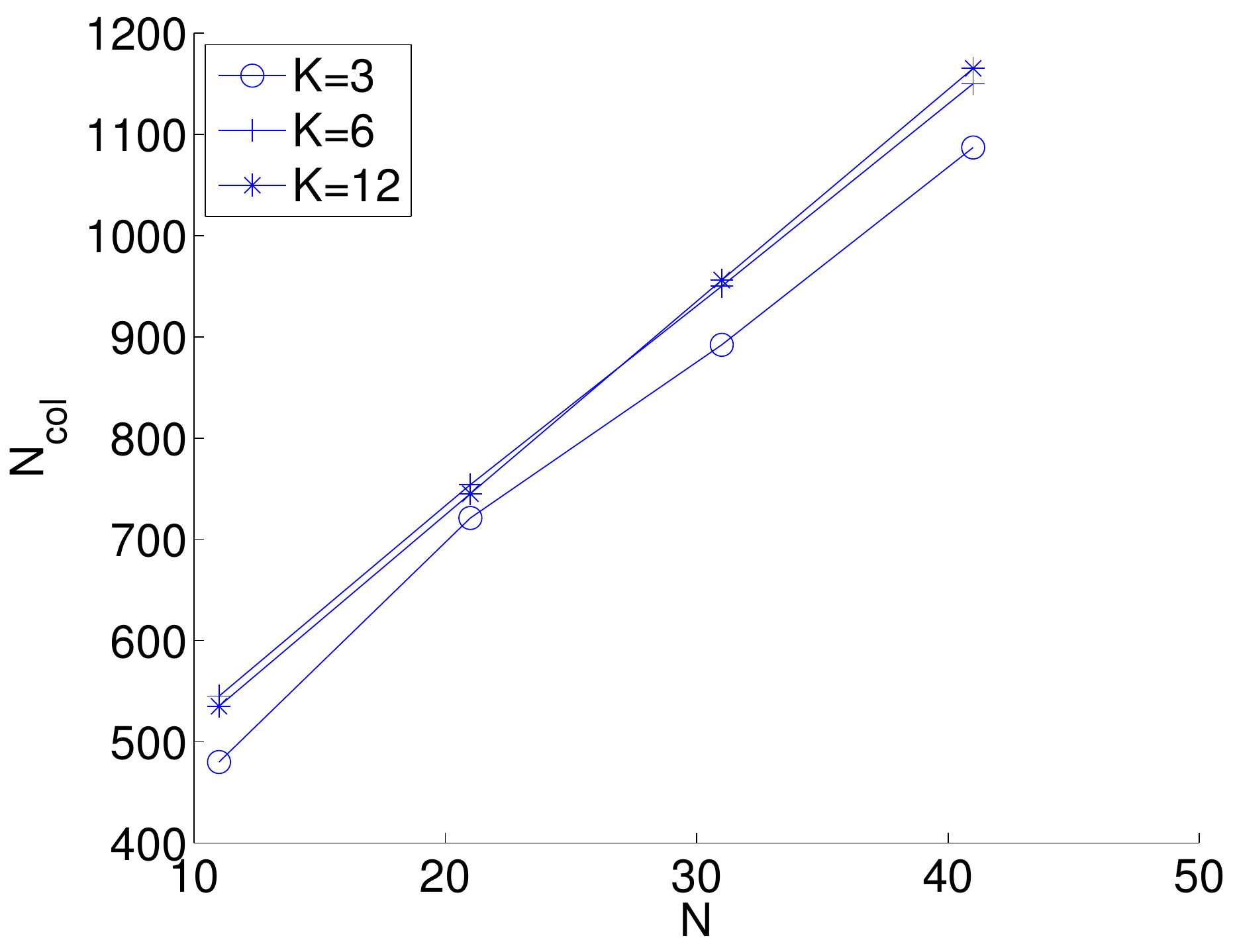} & \includegraphics[height=1.5in]{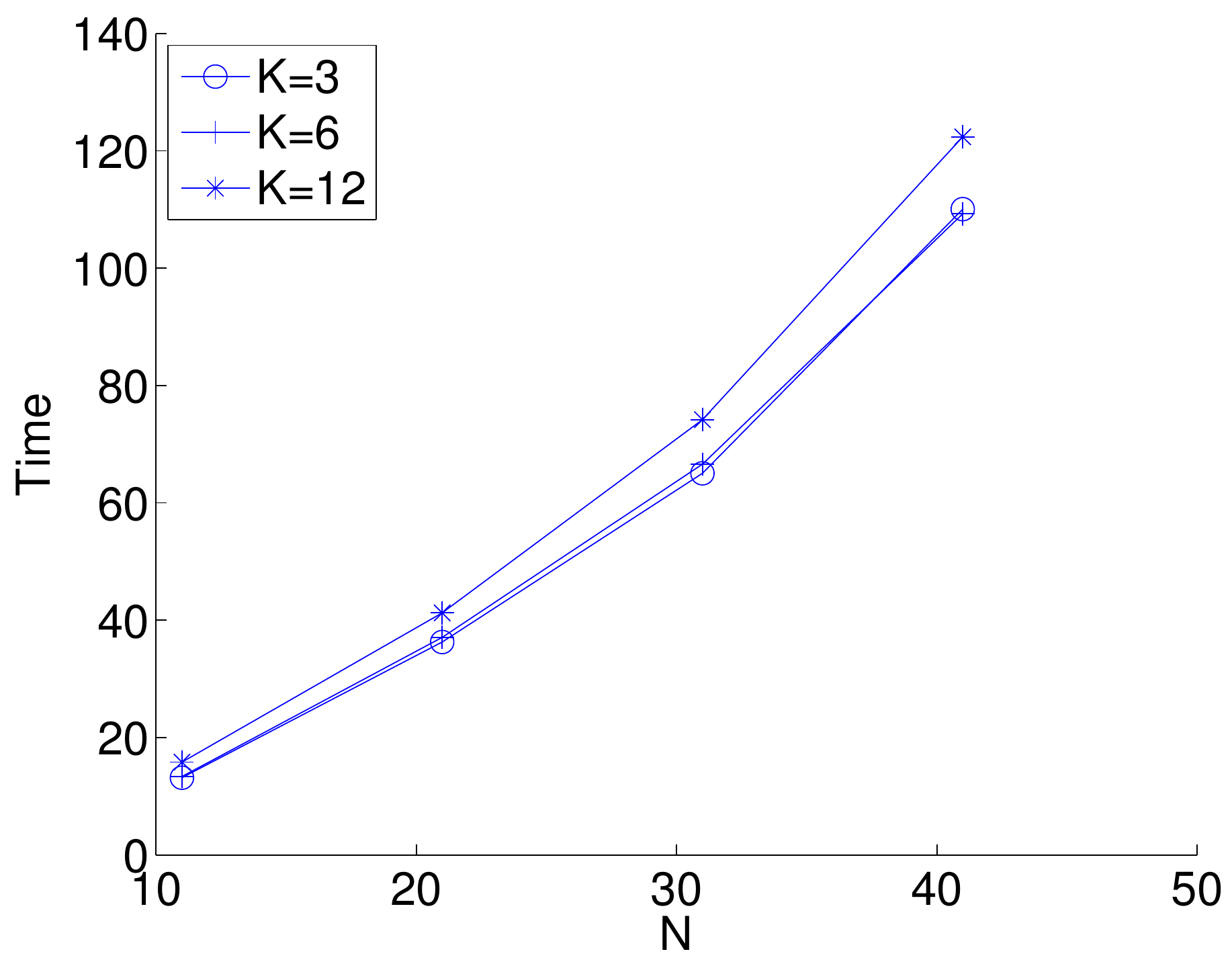} \\
    (c) & (d) \\
    \includegraphics[height=1.5in]{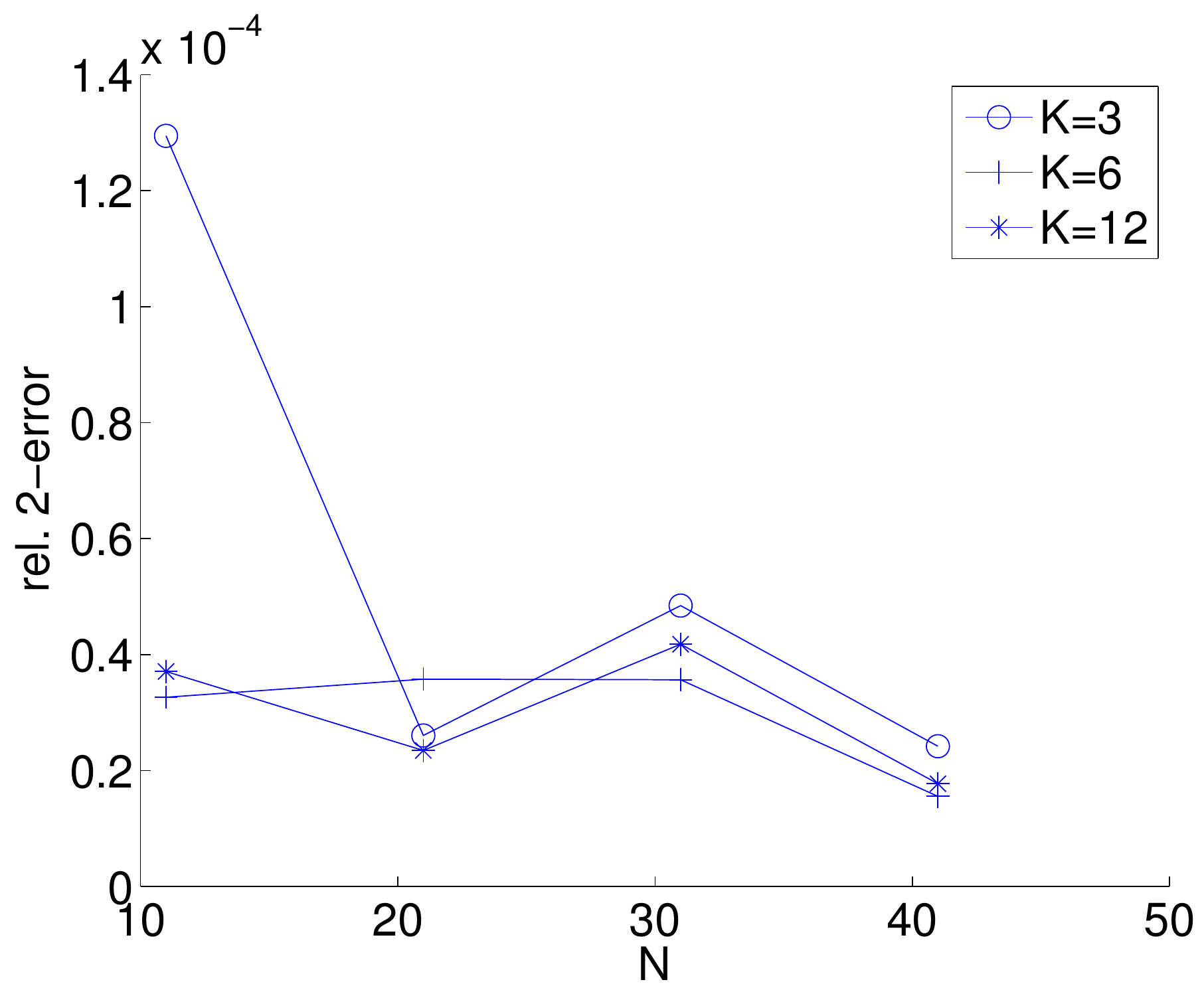} & \includegraphics[height=1.5in]{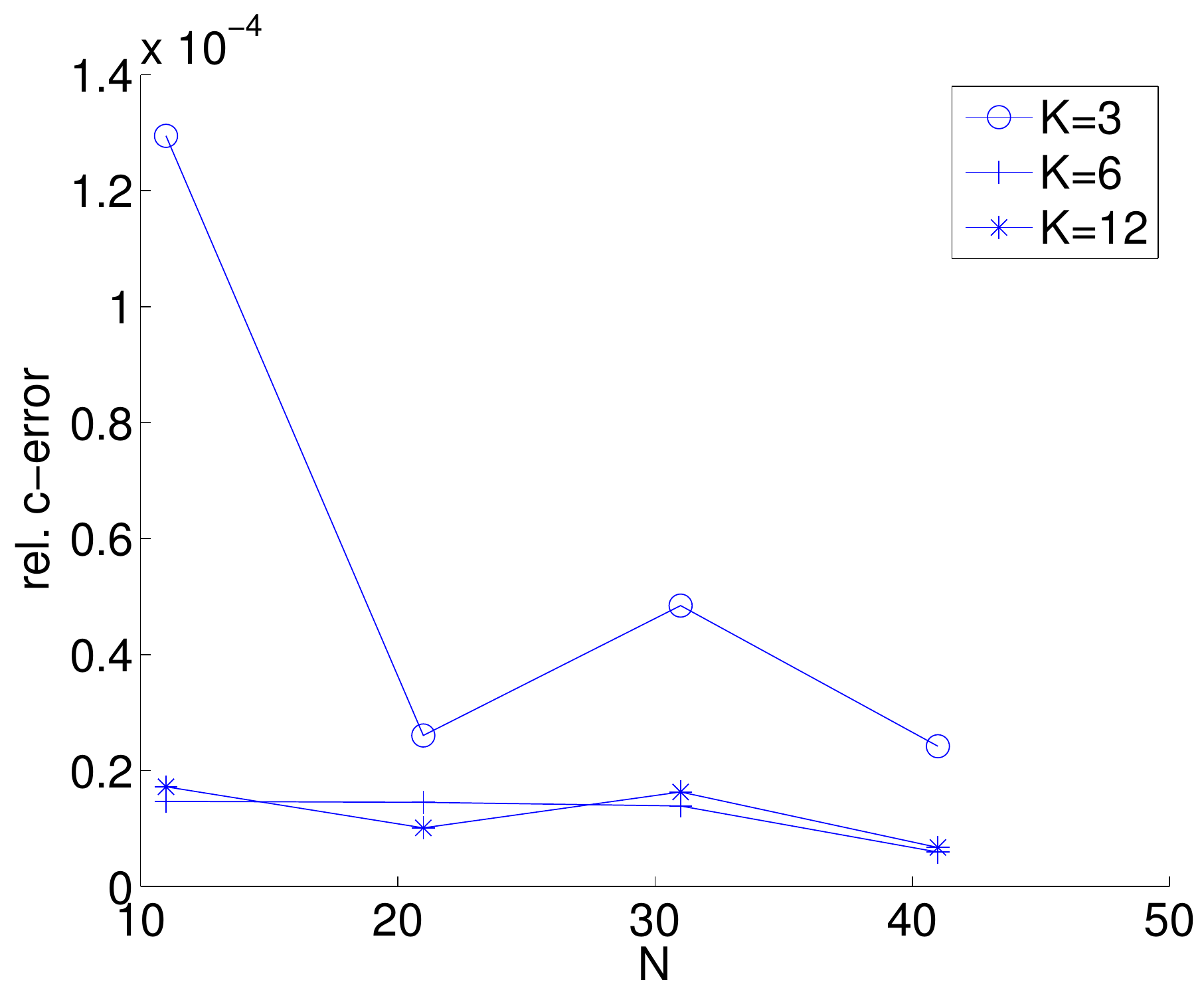}\\
    (e) & (f)
  \end{tabular}
  \caption{Results of the first 3D example.
    (a) the cross-section of the potential $V_c$ in unit cell taken at $x_3=0$.
    (b) the band structure plotted along the $\Gamma-X-M-\Gamma-R-X | M-R$ path.
    (c) $N_{\col}$ as a function $N$ (the number of bands) for different values of $K$ (the number of $k$ points).
    (d) The time used by the column selection algorithm in seconds.
    (e) The relative error measured in the $L^2$ metric.
    (f) The relative error measured in the Coulomb metric.
  }
  \label{fig:t3d4}
\end{figure}

The potential $V_c(x)$ in the first example is given by
\eqref{eq:t2d4} but now in 3D and with $\sigma=0.1667$. Figure
\ref{fig:t3d4} summarizes the results for this examples. The meaning
of each plot is similar to the ones of the 2D examples, except that
\begin{itemize}
\item Figure \ref{fig:t3d4} (a) is the cross-section of the potential
  $V_c(x)$ taken at the plane $x_3=0$;
\item Figure \ref{fig:t3d4} (b) gives the band structure along the
  standard $\Gamma-X-M-\Gamma-R-X | M-R$ path in the Brillouin zone
  for the simple cubic lattice.
\end{itemize}
Similar to the results from the
  2D examples, we can draw the following conclusions:
\begin{itemize}
\item $N_{\col}$ depends roughly linearly on $N$ but is almost
  independent of $K$;
\item the time of the column selection algorithm grows quadratically
  in $N$ but is almost independent of $K$; 
\item the relative errors of the algorithm are bounded by a small
  factor times the accuracy level used in the column selection
  algorithm both in the $L^2$ metric and the Coulomb metric.
\end{itemize}

\begin{figure}[ht!]
  \begin{tabular}{cc}
    \includegraphics[height=1.5in]{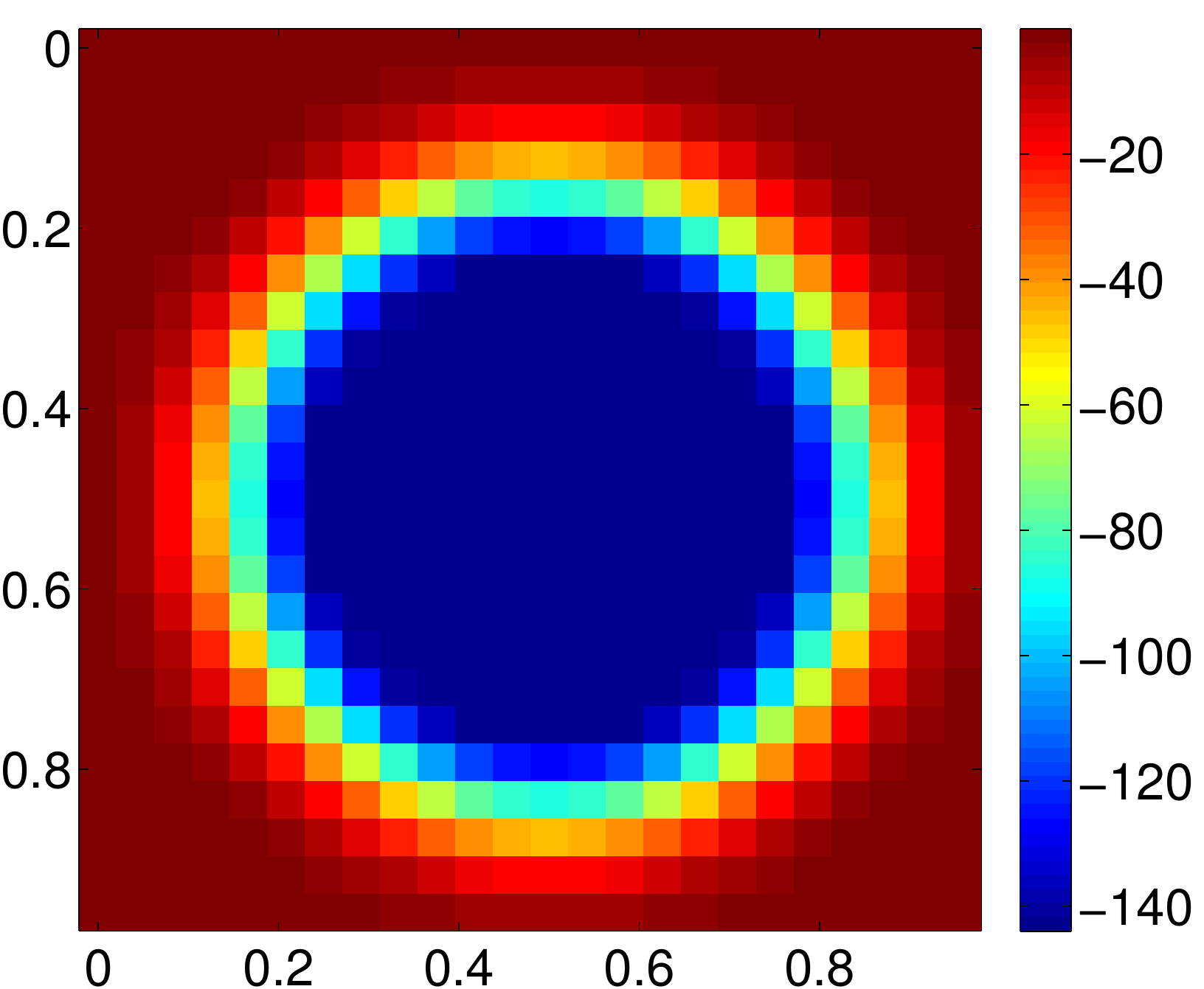} & \includegraphics[height=1.5in]{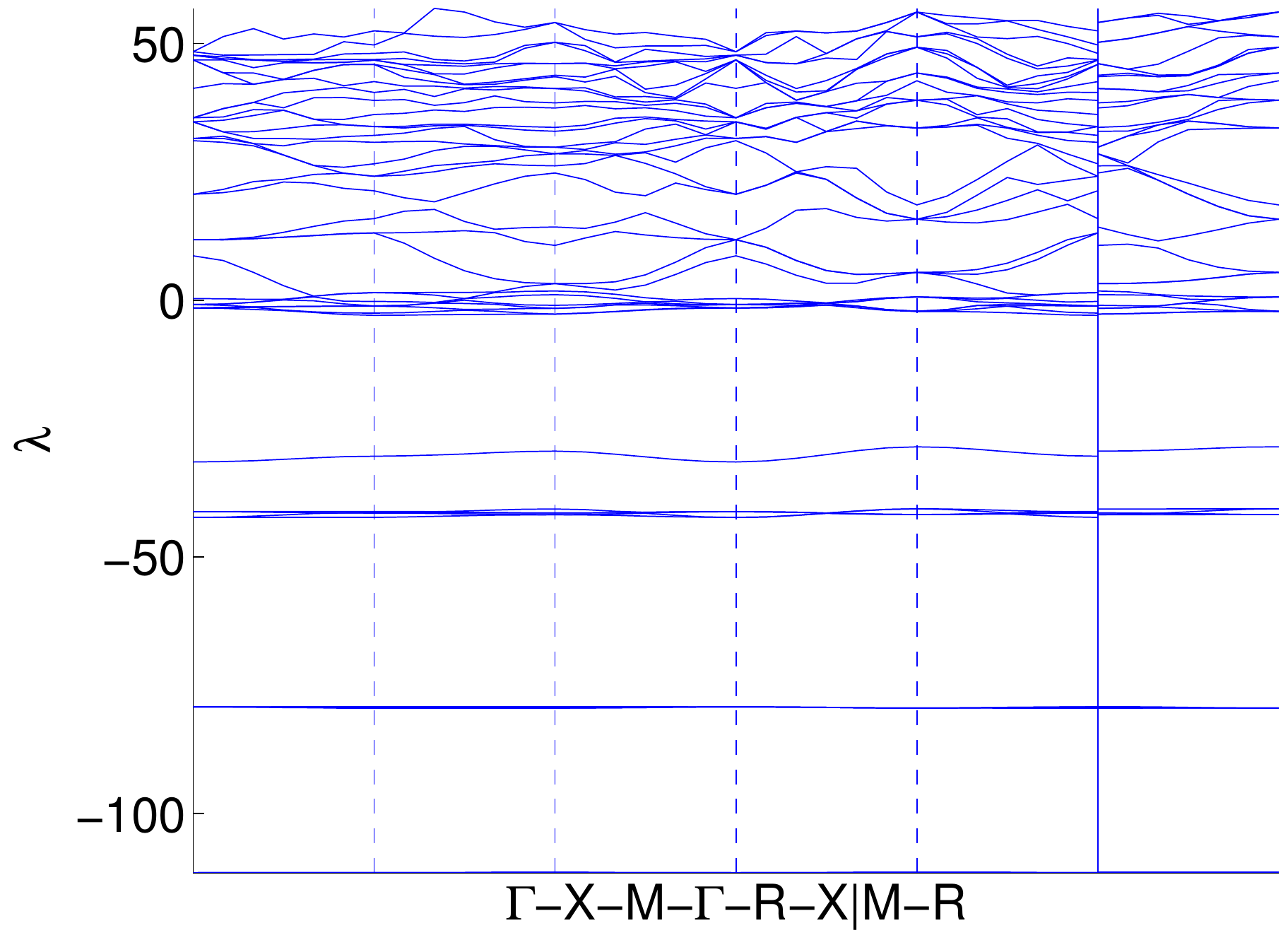} \\
    (a) & (b) \\
    \includegraphics[height=1.5in]{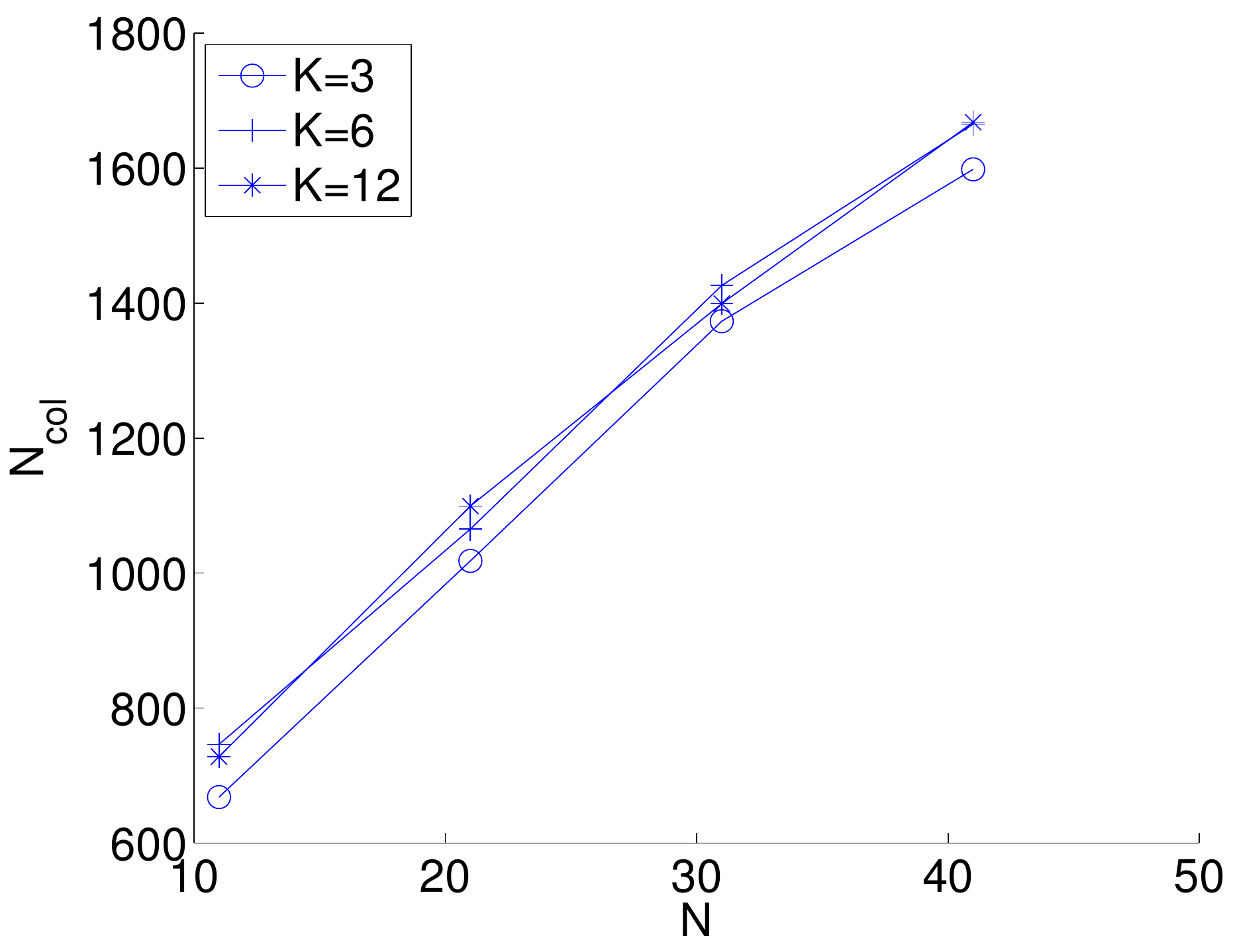} & \includegraphics[height=1.5in]{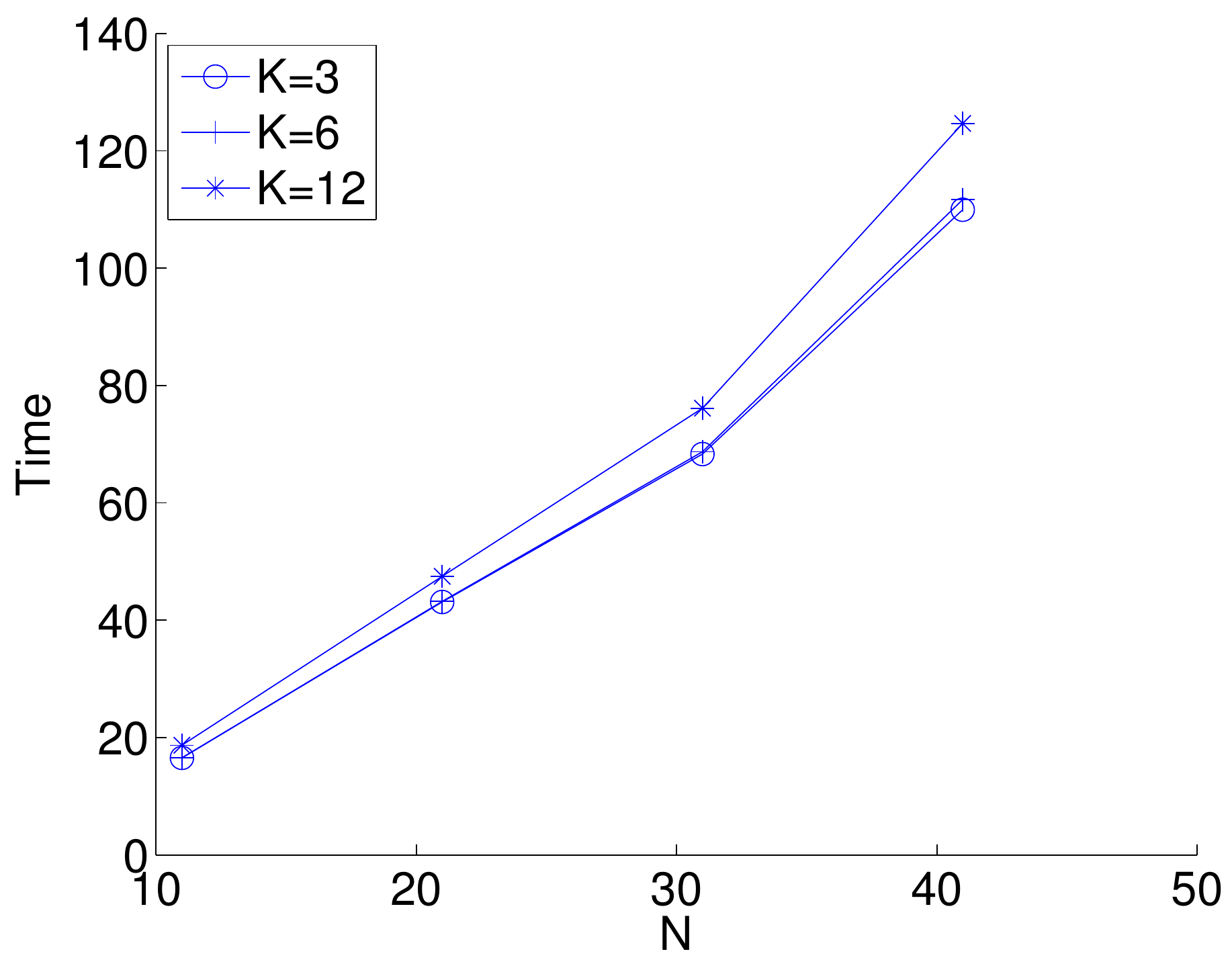} \\
    (c) & (d) \\
    \includegraphics[height=1.5in]{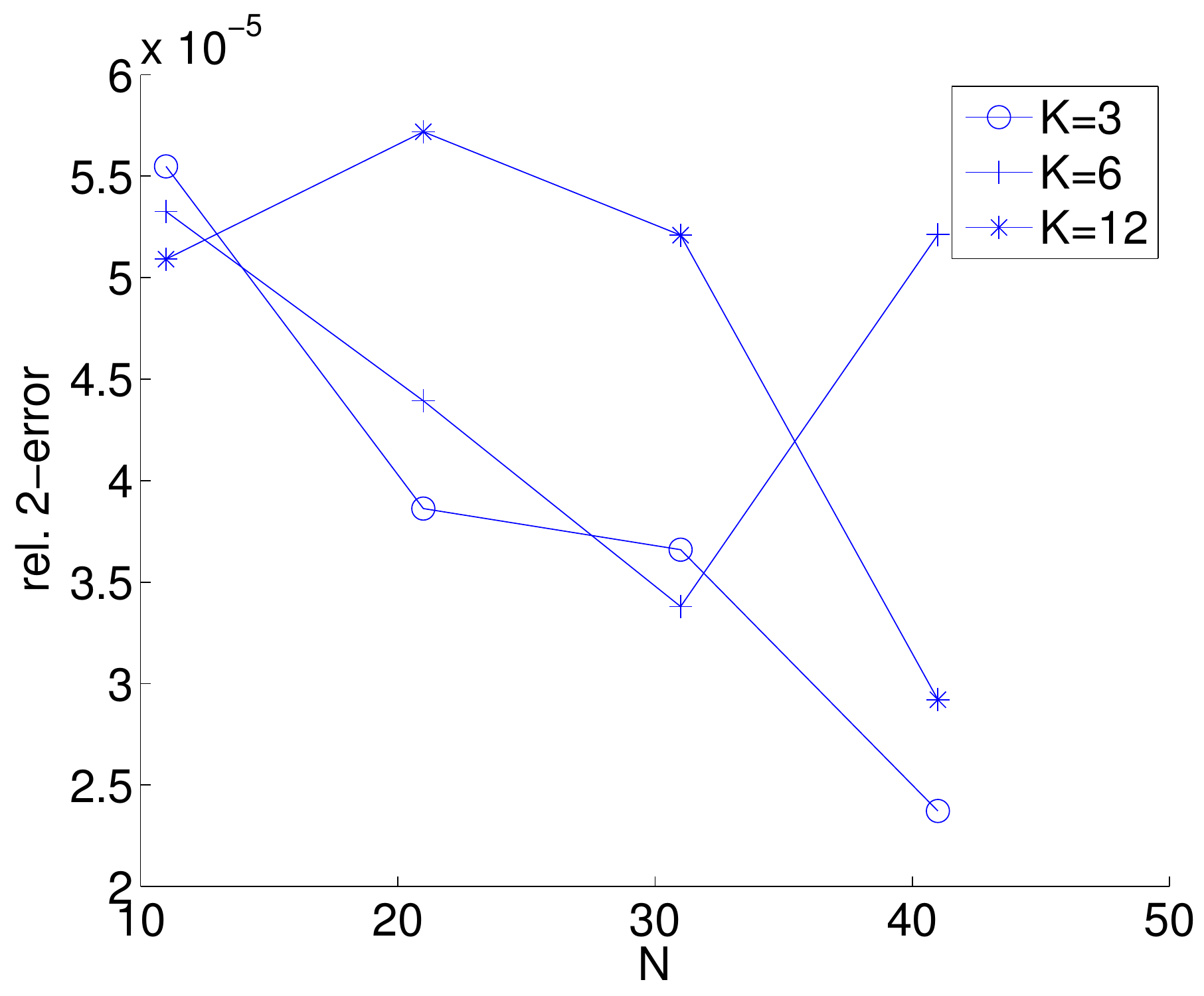} & \includegraphics[height=1.5in]{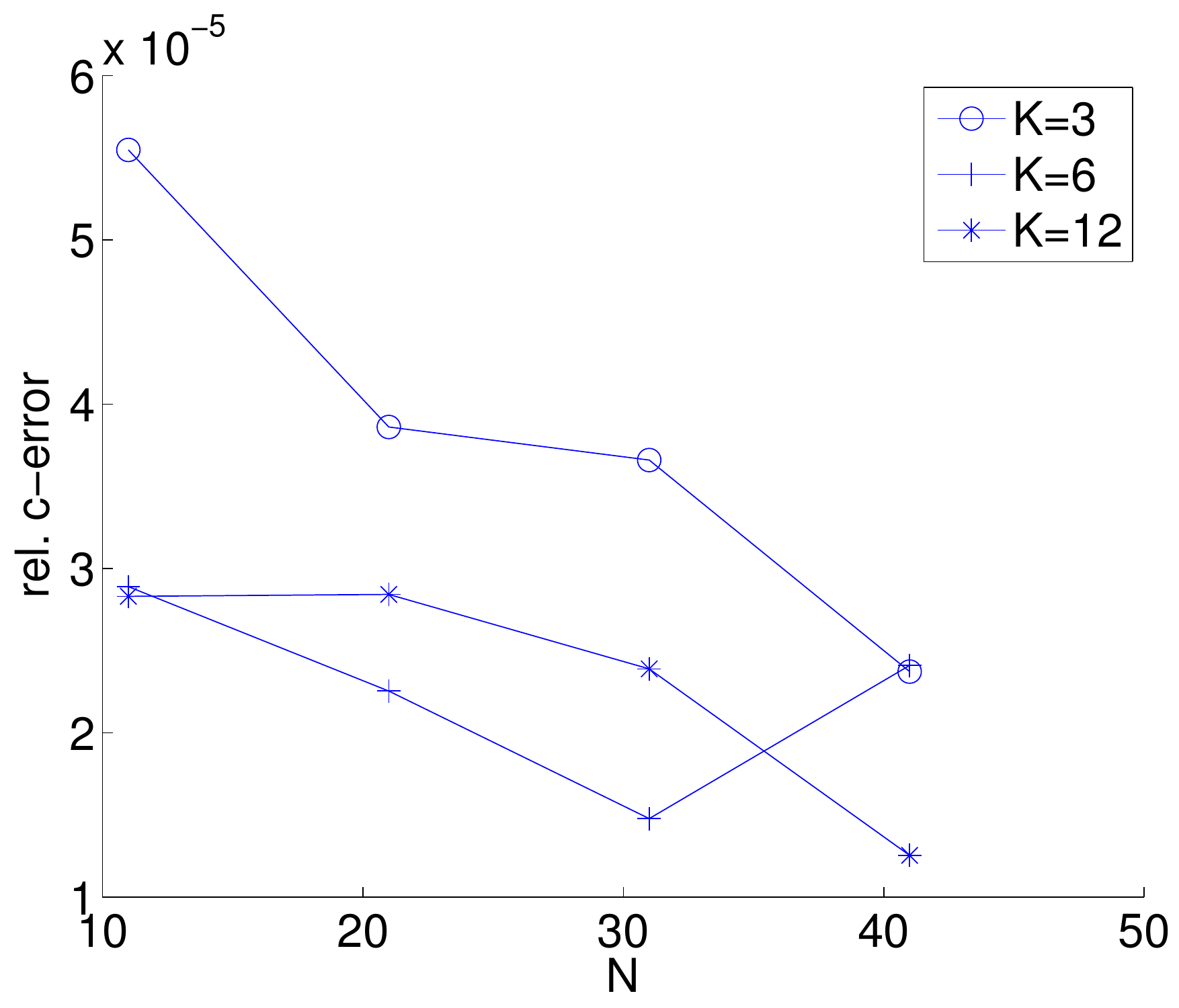}\\
    (e) & (f)
  \end{tabular}
  \caption{Results of the second 3D example.
    (a) the cross-section of the potential $V_c$ in unit cell taken at $x_3=0$.
    (b) the band structure plotted along the $\Gamma-X-M-\Gamma-R-X | M-R$ path.
    (c) $N_{\col}$ as a function $N$ (the number of bands) for different values of $K$ (the number of $k$ points).
    (d) The time used by the column selection algorithm in seconds.
    (e) The relative error measured in the $L^2$ metric.
    (f) The relative error measured in the Coulomb metric.
  }
  \label{fig:t3d5}
\end{figure}

In the second example, the potential $V_c(x)$ given by \eqref{eq:t2d5}
but again interpreted in 3D and with $\sigma=0.0833$. The results are
summarized in Figure \ref{fig:t3d5}.

\section{Conclusion}

We present and demonstrate an efficient algorithm for periodic density
fitting for Bloch waves. The proposed algorithm is based on randomized
pivoted QR algorithm with a choice of random projection adapted to the
tensor product structure of the matrix. The resulting algorithm is
validated in several numerical examples in two and three dimensions.

\bibliographystyle{amsxport}
\bibliography{dft}

\end{document}